\documentclass[reqno]{amsart}
\usepackage{amsmath,amssymb,amsthm,mathtools}
\usepackage{microtype}
\usepackage{tikz}
\usepackage[hidelinks]{hyperref}
\numberwithin{equation}{section}
\newtheorem{theorem}{Theorem}[section]
\newtheorem{proposition}[theorem]{Proposition}
\newtheorem{lemma}[theorem]{Lemma}

\theoremstyle{definition}

\theoremstyle{remark}
\newtheorem{remark}{Remark}[section]
\newcommand{\R}{\mathbb R}
\newcommand{\Sone}{\mathbb S^1}
\newcommand{\A}{\mathcal A}
\newcommand{\Lop}{\mathcal L}
\newcommand{\av}[1]{\overline{#1}}
\newcommand{\dd}{\,d}

\newcommand{\Span}{\operatorname{span}}
\newcommand{\Ran}{\operatorname{Ran}}
\newcommand{\Ker}{\operatorname{Ker}}

\newcommand{\Id}{\operatorname{Id}}

\hypersetup{
 pdftitle={Bifurcation for annular overdetermined problems with linear and power nonlinearities},
 pdfauthor={Jun Wang, Yong Zhang},
 pdfkeywords={Annular domains; overdetermined elliptic problems; power nonlinearities; bifurcation; thin annuli}
}

\begin{document}
\raggedbottom
\title[Annular overdetermined problems with power nonlinearities]
{On the annular overdetermined
problems with linear and power nonlinearities in $\mathbb{R}^2$}

\author[J. Wang]{Jun Wang}
\address{Jun Wang, School of Mathematical Sciences, Jiangsu University,
Zhenjiang, Jiangsu 212013, P.R. China}
\email{wangmath2011@126.com}
\author[F. Xu]{Fei Xu}
\address{Fei Xu, School of Mathematical Sciences, Jiangsu University,
Zhenjiang, Jiangsu 212013, P.R. China}
\email{xufeiujs@126.com}
\author[Y. Zhang]{Yong Zhang}
\address{Yong Zhang, School of Mathematical Sciences, Jiangsu University,
Zhenjiang, Jiangsu 212013, P.R. China}
\email{zhangyong@ujs.edu.cn}

\begin{abstract}
In this paper, we study two classes of partially overdetermined elliptic
problems on planar annular domains in $\mathbb{R}^2$, with a fixed inner circle and a free
outer boundary. The Dirichlet values are prescribed to be $1$ and $0$
on the inner and outer boundaries, respectively, while the squared
normal derivative is constant on the outer boundary. For the linear law
$\Delta u=-\mu u$, we establish local bifurcation under a Dirichlet
nonresonance condition. For the power law $\Delta u=u-u^p$ with $p>1$,
we take the mean outer radius as the bifurcation parameter and construct
large positive radial solutions on thin annuli. The limiting profile
solves $-V''=V^p$ on an interval, and the principal mixed
Dirichlet--Neumann eigenvalue of its linearization determines a constant
$\kappa_p>0$. For every fixed inner radius $\lambda>0$ and every
sufficiently large integer $K$, we obtain a critical radius
$R_K=\lambda+\lambda\kappa_p/K+o(K^{-1})$ and a smooth local curve
of positive solutions with noncircular, even, $2\pi/K$-periodic outer
boundaries. The proof combines shape derivatives, spectral estimates,
and the Crandall--Rabinowitz bifurcation theorem.
\end{abstract}
\maketitle

\section{Introduction and main results}

In this paper, we study partially overdetermined elliptic problems on annular domains with a fixed inner circle and a free outer boundary. We are interested in the parameter values at which nontrivial domains bifurcate from concentric annuli. The classical starting point is Serrin's overdetermined problem, in which a positive solution has zero Dirichlet data and constant normal derivative on the boundary of a bounded domain. Serrin's moving-plane argument \cite{Serrin} and Weinberger's proof \cite{Weinberger} show that the domain must be a ball. For an annular domain whose boundary components carry distinct Dirichlet values, additional boundary data may instead permit nonradial solutions. In \cite{GWYZ}, the authors investigated this possibility through the spectrum of the linearized boundary operator for a constant vorticity function. Related annular bifurcation results were obtained by Kamburov and Sciaraffia \cite{KS}, who also considered constant vorticity and imposed overdetermined data on both boundary components. The present paper extends this line of research to a nonconstant vorticity function.

We first mention some related work in this area. For unbounded domains, Schlenk and Sicbaldi \cite{SS} constructed a smooth family of noncompact domains bifurcating from the straight cylinder $B^n\times\mathbb{R}$. On these domains, they found only a positive solution of the linear equation $\Delta v+\sigma v=0$ with zero Dirichlet data and constant normal derivative on the boundary, while Dai and Zhang \cite{DZ} found sign-changing solutions. These perturbed domains constructed above are rotationally symmetric and periodic in the axial direction. Ruiz, Sicbaldi and Wu \cite{RSW} extended this cylindrical bifurcation approach to semilinear equations $\Delta v+g(v)=0$ with general nonlinearities. Starting from a positive radial Dirichlet solution in the cross-sectional ball, and under their nondegeneracy hypotheses, they constructed positive bounded solutions on nontrivial periodic, onduloid-type domains bifurcating from a straight cylinder. Fall, Minlend and Weth \cite{FMW} studied the Schiffer problem ($\Delta v+\sigma v=0$ with zero Neumman data and constant Dirichlet data on the boundary) on the flat cylinder $\mathbb{R}^N\times(\mathbb{R}/2\pi\mathbb{Z})$ and on the sphere $\mathbb{S}^2$. They constructed domains carrying Neumann eigenfunctions with a single constant boundary trace and whose boundaries have nonconstant principal curvatures; their spherical examples disprove Souam's conjecture. Their functional framework overcomes the loss of regularity that arises in standard spaces.

For bounded domains, Ruiz \cite{RuizJEMS} showed that the positivity assumption in Serrin's symmetry theorem cannot simply be omitted. For a suitable $C^2$ nonlinearity, he constructed a sign-changing solution of $-\Delta v=v-\left(v^+\right)^3$ on a smooth bounded domain that is not a ball, with zero Dirichlet data and constant normal derivative on the whole boundary. The construction uses local bifurcation and the spectrum of the associated linearized operator. Thus the sign of the solution is an essential part of the symmetry question. Recent work on Schiffer-type problems also illustrates how the boundary conditions and the ambient geometry affect symmetry. Enciso, Fern\'andez, Ruiz and Sicbaldi \cite{EFRS} constructed noncircular doubly connected planar domains bifurcating from annuli that support a nonradial solution of $\Delta v+\mu v=0$ with zero normal derivative and constant trace on each boundary component, the two trace values being allowed to differ. Both boundary components are perturbed, and the construction also yields continuous, compactly supported stationary weak Euler flows that are not locally radial.
For another different Helmholtz overdetermined problem, Wheeler \cite{Wheeler} constructed real-analytic branches of noncircular bounded planar domains near a disk, with $v=b$ and $\partial_\nu v=c$ on the whole boundary, where $b$ and $c$ are nonzero constants. For each symmetry order $m\ge4$, his argument combines conformal coordinates, the Crandall--Rabinowitz theorem and a simple zero of a Bessel-function Wronskian, providing counterexamples to the Willms--Gladwell conjecture.

These constructions concern, respectively, vanishing Neumann data on annular or curved domains and nonzero constant Cauchy data on a simply connected planar domain.They lead naturally to the question considered here: whether bifurcation can occur for a Bernoulli-type overdetermined elliptic problem, that is to say, whether a domain with a noncircular outer boundary can bifurcate from a concentric annulus when the inner circle is fixed, the two Dirichlet values are prescribed separately, and nonzero Neumann data is imposed only on the outer boundary.
Here we will answer this equation combing shape derivatives, spectral estimates,
and the Crandall--Rabinowitz bifurcation theorem.
More precisely, let $\lambda>0$ and let $\Omega$ be a bounded domain
containing $\overline{B_\lambda}$, where $B_{\lambda}$ is a disk of radius
$\lambda$. Set
$\A_\Omega=\Omega\setminus\overline{B_\lambda}$. We consider the
system
\begin{equation}\label{eq:general}
 \begin{cases}
  \Delta u=f(u)&\text{in }\A_\Omega,\\
  u=1&\text{on }\partial B_\lambda,\\
  u=0&\text{on }\partial\Omega,\\
  |\nabla u|^2=Q>0&\text{on }\partial\Omega.
 \end{cases}
\end{equation}
Throughout the analysis, the inner circle remains fixed; that is, $\lambda$ is prescribed. We describe
the perturbed outer boundary by
\begin{equation}\label{eq:geometry}
 \partial\Omega_{R,\eta}
 =\{(R+\eta(\theta))(\cos\theta,\sin\theta):\theta\in\Sone\},
 \qquad \av\eta:=\frac1{2\pi}\int_0^{2\pi}\eta(\theta)\dd\theta=0.
\end{equation}
Here $R$ denotes the mean radial coordinate of the outer boundary.
The zero-mean condition fixes the decomposition into the radius and
the boundary perturbation.
We begin with the linear-vorticity problem
\begin{equation}\label{eq:linear}
 \begin{cases}
  \Delta u=-\mu u&\text{in }\Omega\setminus\overline{B_\lambda},\\
  u=0&\text{on }\partial\Omega,\\
  u=1&\text{on }\partial B_\lambda,\\
  |\nabla u|^2=Q>0&\text{on }\partial\Omega.
 \end{cases}
\end{equation}
For this problem, the reference outer radius is $1$, and $\mu$ is
the bifurcation parameter.
We compare \eqref{eq:linear} with the Schiffer-type problem studied by
Enciso, Fern\'andez, Ruiz and Sicbaldi \cite{EFRS}. They construct
noncircular doubly connected domains carrying a nonradial solution of
$\Delta v+\mu v=0$, with $\partial_\nu v=0$ on both boundary
components and a constant trace on each component, the two constants
being allowed to differ. Their construction perturbs both circles of
a reference annulus and yields continuous, compactly supported
stationary weak Euler flows that are not locally radial. In our
problem, we keep the inner circle fixed, prescribe the boundary values
$1$ and $0$, and impose $|\partial_\nu u|=\sqrt Q>0$ only on the
outer boundary. We also use bifurcation from concentric annuli, but
the relevant spectral conditions are different: their argument relies
on a crossing between a radial Neumann eigenvalue and a nonradial
Dirichlet eigenvalue, whereas our linearized boundary condition leads
to a mixed Dirichlet--Robin eigenvalue problem.

We then consider the power-law problem for
\begin{equation}\label{eq:power}
 \begin{cases}
  \Delta u=u-u^p&\text{in }\Omega\setminus\overline{B_\lambda},\\
  u>0&\text{in }\Omega\setminus\overline{B_\lambda},\\
  u=0&\text{on }\partial\Omega,\\
  u=1&\text{on }\partial B_\lambda,\\
  |\nabla u|^2=Q>0&\text{on }\partial\Omega,
 \end{cases}
\end{equation}
where $p>1$. For this problem, the reference outer radius $R$
is undetermined and serves as the bifurcation parameter. The main equation $\Delta u=u-u^p$, equivalently
$-\Delta u+u-u^p=0$, also appears in the exterior-domain construction
of Ros, Ruiz and Sicbaldi \cite{RRS}.
Ros, Ruiz and Sicbaldi construct exterior domains admitting positive
solutions with zero Dirichlet data and constant normal derivative.
Their bifurcation parameter is a radius, which becomes a diffusion
coefficient after rescaling, and their spectral analysis is based on
the positive radial exterior solution. In \eqref{eq:power}, the domain
is bounded and the inner boundary has the nonzero Dirichlet value $1$.
These boundary conditions require a separate construction of the radial
solution and a separate verification of the spectral gap. We establish
both properties, together with the boundary eigenvalue crossing, by
analyzing the thin-annulus limit.

For the local analysis, we use the standard H\"older spaces
$C^{m,\alpha}(\Sone)$, where $m\ge0$ is an integer and
$0<\alpha<1$, endowed with their usual H\"older norms, and the
corresponding spaces on smooth compact domain closures.
For each integer $K\ge1$, define
\[
 X_K=\{\eta\in C^{2,\alpha}(\Sone):
 \eta(-\theta)=\eta(\theta),\
 \eta(\theta+2\pi/K)=\eta(\theta),\ \av\eta=0\},
\]
and let
\[
 Y_K=\{\eta\in C^{1,\alpha}(\Sone):
 \eta(-\theta)=\eta(\theta),\
 \eta(\theta+2\pi/K)=\eta(\theta),\ \av\eta=0\}.
\]
The nonzero Fourier modes in these spaces are
$\cos(nK\theta)$, $n\ge1$. We first state the main result for the
linear problem. For $\lambda\in(0,1)$ and $K\ge1$, we denote by
$\mu_{K,\lambda}^*$ the lowest eigenvalue of
\begin{equation}\label{eq:linear-Robin}
 -h''-\frac1r h'+\frac{K^2}{r^2}h=\mu h,\qquad
 h(\lambda)=0,\quad h'(1)+h(1)=0
\end{equation}
and denote by $\Sigma_{\mathrm{rad}}(\lambda)$ the radial Dirichlet spectrum of $-\Delta$ on $B_1\setminus\overline{B_\lambda}$.
We first state the annular bifurcation result for the linear vorticity law.

\begin{theorem}\label{thm:linear}
Fix $\lambda\in(0,1)$ and $K\ge1$, and suppose that
$\mu_{K,\lambda}^*\notin\Sigma_{\mathrm{rad}}(\lambda)$.
Then there exist $s_0>0$ and a $C^\infty$ curve
\[
 s\longmapsto(\mu(s),\eta(s),Q(s)),\qquad |s|<s_0,
\]
of solutions to \eqref{eq:linear} on
$\Omega_{1,\eta(s)}\setminus\overline{B_\lambda}$, with
$\eta(s)\in X_K$, satisfying
\[
 \mu(s)=\mu_{K,\lambda}^*+O(s^2),\qquad
 \eta(s)=s\cos(K\theta)+O(s^2),\qquad Q(s)>0.
\]
The inner boundary remains $\partial B_\lambda$ along the curve.
\end{theorem}
\begin{remark}
The nonresonance condition $\mu_{K,\lambda}^*\notin\Sigma_{\mathrm{rad}}(\lambda)$ is automatic for $K=1$. For each fixed $K\ge2$, it fails only on a discrete set of inner radii in $(0,1)$ and holds for all $\lambda$ sufficiently close to $1$. Moreover, in the chosen symmetry class, only the radial mode can cause Dirichlet resonance, as proved in Lemma~\ref{lem:linear-nonresonance}. Thus, for $K=1$, we obtain a nontrivial deformation of the outer boundary while keeping the inner circle fixed without imposing any additional assumption. In addition, we would like to mention that the fixed inner circle breaks translation invariance, and the kernel $\cos\theta$ comes from the mixed eigenvalue problem \eqref{eq:linear-Robin}.
\end{remark}

\begin{remark}[Rigidity for the linear problem]\label{rem:linear-rigidity}
Let $\Omega$ be a bounded simply connected $C^2$ domain containing
$\overline{B_\lambda}$, and let $u$ be a classical solution of
\eqref{eq:linear}. If $0<u<1$ in $\A_\Omega$, then
$\Omega=B_R(0)$ for some $R>\lambda$ and $u(x)=U(|x|)$, with
$U'(r)<0$ for $\lambda<r<R$. On the inner circle, $u=1$ and
$u<1$ in $\A_\Omega$ imply that the radial derivative pointing
into the annulus is nonpositive. The symmetry conclusion follows from
Reichel's symmetry theorem \cite[Theorem~3, case~(15b)]{Reichel},
since the outer boundary condition gives
$\partial_\nu u=-\sqrt Q$ on $\partial\Omega$. In particular,
when $\mu\le0$, the maximum principle gives $0<u<1$ automatically,
so every classical solution is radial on a concentric annulus.
Under the mean outer-radius normalization used in
Theorem~\ref{thm:linear}, we have $R=1$.
No assumption that the domain is close to an annulus is needed.
\end{remark}

For the power-law problem \eqref{eq:power}, we assume
\begin{equation}\label{eq:alpha}
 0<\alpha<\min\{1,p-1\}.
\end{equation}
This restriction gives the H\"older regularity needed for nonintegral
powers at the outer boundary.
To state the main result on the problem \eqref{eq:power}, we first introduce the constant
governing the critical radii.
Let $V=V_p$ be the unique positive solution of
\[
 -V''=V^p\quad\text{in }(0,1),\qquad V(0)=V(1)=0.
\]
We prove its existence and uniqueness in Section~\ref{sec:pulse}.
Denote by $\nu_{\mathrm{DN}}(p)$ the lowest eigenvalue of
\begin{equation}\label{eq:DN-intro}
 -h''-pV^{p-1}h=\nu h,\qquad h(0)=0,\quad h'(1)=0,
\end{equation}
and set
\begin{equation}\label{eq:kappa}
 \kappa_p=\sqrt{-\nu_{\mathrm{DN}}(p)}.
\end{equation}
We shall prove in Section~\ref{sec:pulse} that
$\nu_{\mathrm{DN}}(p)<0$. Thus $\kappa_p$ is well defined.
The main result for problem \eqref{eq:power} is as follows.
We next state the annular bifurcation result for the power law.

\begin{theorem}\label{thm:power}
Fix $\lambda>0$, $p>1$, and $\alpha$ satisfying \eqref{eq:alpha},
there exists an integer $K_0=K_0(\lambda,p)\ge2$ such that, for
each integer $K\ge K_0$, there exist a radius $R_K>\lambda$,
a number $s_K>0$, and a $C^\infty$ curve
\[
 s\longmapsto (R(s),\eta(s),Q(s)),\qquad |s|<s_K,
 \qquad \eta(s)\in X_K,
\]
of solutions to \eqref{eq:power} on
$\Omega_{R(s),\eta(s)}\setminus\overline{B_\lambda}$ satisfying
\[
 R(s)=R_K+O(s^2),\qquad
 \eta(s)=s\cos(K\theta)+O(s^2),\qquad Q(s)>0.
\]
Moreover, the critical radii satisfy
\begin{equation}\label{eq:RK-asymptotic}
 R_K=\lambda+\frac{\lambda\kappa_p}{K}+o(K^{-1}).
\end{equation}
Let $U_K$ denote the radial state at $s=0$. There is a neighborhood
of $(R_K,0,U_K)$, measured in the outer-radius variable, the boundary
space $X_K$, and the $C^{2,\alpha}$ norm of the state pulled back
to the fixed annulus $B_{R_K}\setminus\overline{B_\lambda}$ by
\eqref{eq:power-pullback} with $R_*=R_K$, such that every nonradial
solution in this neighborhood, with boundary perturbation lying in
$X_K$, belongs to this local curve.
\end{theorem}

\begin{remark}[Rigidity for the power-law problem]\label{rem:power-rigidity}
Under the geometric assumptions of Remark~\ref{rem:linear-rigidity},
any classical solution of \eqref{eq:power} satisfying $u\le1$ in
$\A_\Omega$ has $\Omega=B_R(0)$ for some $R>\lambda$ and
$u(x)=U(|x|)$, with $U'<0$ on $(\lambda,R)$. Indeed, the strong
maximum principle gives $0<u<1$. The radial derivative pointing
into the annulus is nonpositive on the inner circle, as in
Remark~\ref{rem:linear-rigidity}, and the conclusion follows from
\cite[Theorem~3, case~(15b)]{Reichel}.
\end{remark}

In Theorem~\ref{thm:power}, the inner radius is fixed along each branch, while
the mean outer radius varies. The proof constructs the critical radii
for this fixed nonlinear equation and verifies the boundary eigenvalue
crossing.
The nonlinear analysis requires two additional observations. First,
the potential in the linearized equation is $1-pU^{p-1}$, where
the radial solution $U$ depends on the parameter. The transversality
calculation must therefore include the derivative of $U$. Second,
the inner boundary value $1$ is an equilibrium of $u-u^p$.
This property excludes the first harmonic from the boundary kernel,
as stated below.

\begin{proposition}\label{prop:first}
Let $U$ be a positive radial solution of
\[
 U''+r^{-1}U'=U-U^p,\qquad U(\lambda)=1,\quad U(R)=0,
 \qquad R>\lambda.
\]
Then the boundary-value problem
\[
 -h''-r^{-1}h'
 +\bigl(r^{-2}+1-pU^{p-1}\bigr)h=0,\qquad
 h(\lambda)=0,\quad h'(R)+R^{-1}h(R)=0
\]
admits only the zero solution. When the corresponding Dirichlet state
map is defined, the first harmonic therefore does not belong to the
kernel of the linearized Bernoulli operator at this radial solution.
\end{proposition}

We prove this proposition at the end of Section~\ref{sec:nonlinear-proof} by using
$U'$ as a Jacobi field. Thus, for the power law, the bifurcation
conclusion cannot extend to every mode $K\ge1$. The sufficiently
large modes in Theorem~\ref{thm:power} are compatible with this
first-harmonic obstruction.

\begin{figure}[ht]
\centering
\begin{tikzpicture}[scale=1.02]
 \fill[gray!13] (0,0) circle (1.55);
 \fill[white] (0,0) circle (0.90);
 \draw[thick] (0,0) circle (0.90);
 \draw[thick] (0,0) circle (1.55);
 \draw[dashed] (0,0)--(1.55,0);
 \node at (0,0) {$B_\lambda$};
 \node at (0,-1.95) {radial annulus};
 \node at (1.17,0.22) {$\delta$};
 \draw[->,thick] (2.0,0)--(3.0,0);
 \begin{scope}[xshift=5cm]
  \fill[gray!13,domain=0:360,samples=181,smooth]
   plot ({(1.55+0.10*cos(6*\x))*cos(\x)},
         {(1.55+0.10*cos(6*\x))*sin(\x)})--cycle;
  \fill[white] (0,0) circle (0.90);
  \draw[thick] (0,0) circle (0.90);
  \draw[thick,domain=0:360,samples=181,smooth]
   plot ({(1.55+0.10*cos(6*\x))*cos(\x)},
         {(1.55+0.10*cos(6*\x))*sin(\x)});
  \node at (0,0) {$B_\lambda$};
  \node at (0,-1.95) {outer-boundary bifurcation};
 \end{scope}
\end{tikzpicture}
\caption{Schematic illustration of the domains in the nonlinear
bifurcation theorem. The inner circle is fixed, and the outer boundary
develops an angular mode.}
\label{fig:geometry}
\end{figure}
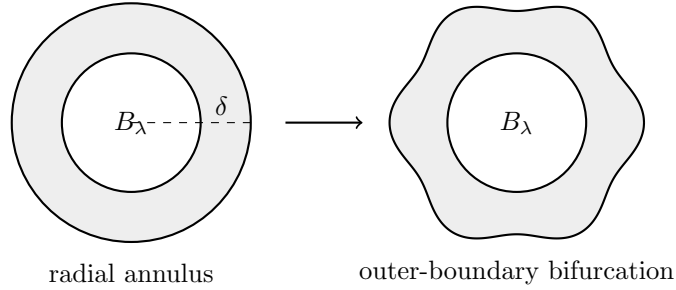

The remainder of the paper is organized as follows.
In Section~\ref{sec:linear}, we prove Theorem~\ref{thm:linear} for
the linear-vorticity problem. We introduce the working spaces and
radial branch, define the state and Bernoulli equations, and compute
the boundary linearization. We then determine its kernel and range
and verify the Crandall--Rabinowitz conditions.
In Section~\ref{sec:nonlinear-proof}, we follow the same approach
for the power-law problem. We construct the large radial solutions
from the one-dimensional limiting profile, define the corresponding
operators, and use the thin-annulus spectral estimates to verify the
bifurcation conditions and prove Theorem~\ref{thm:power}.
We conclude this section with the proof of Proposition~\ref{prop:first}.

\section{The linear-vorticity overdetermined problem}\label{sec:linear}

In this section, we fix the inner radius $\lambda\in(0,1)$
and the mean outer radius $1$, construct the Dirichlet solution for
each small boundary graph, and use its squared outer normal derivative
as the bifurcation equation. Now, let us introduce the working spaces and
fundamental domain.

\subsection{Working spaces and the fixed annulus}\label{sec:fixed}
Fix $0<\alpha<1$ and an integer $K\ge1$, and put
$A_\lambda=B_1\setminus\overline{B_\lambda}$. We write $D_K$ for
the group generated by $\theta\mapsto-\theta$ and
$\theta\mapsto\theta+2\pi/K$. For functions on the annulus,
$C^{m,\alpha}_{D_K}(\overline A_\lambda)$ denotes the closed
subspace invariant under these two transformations. The boundary
spaces used in Theorem~\ref{thm:linear} are
\begin{equation}\label{eq:linear-spaces}
 \begin{split}
 X_K&=\{\eta\in C^{2,\alpha}_{D_K}(\Sone):\av\eta=0\},\\
 Y_K&=\{g\in C^{1,\alpha}_{D_K}(\Sone):\av g=0\}.
 \end{split}
\end{equation}
For the interior equation, we use
\[
 Z_K^{\mathrm{lin}}
 =\{w\in C^{2,\alpha}_{D_K}(\overline A_\lambda):
       w|_{r=\lambda}=w|_{r=1}=0\},\qquad
 E_K^{\mathrm{lin}}=C^{0,\alpha}_{D_K}(\overline A_\lambda).
\]
In particular, $Z_K^{\mathrm{lin}}$ contains radial functions. Only
the boundary graph has zero angular mean; we do not remove the radial
mode from the state equation.

Choose $\chi\in C^\infty([\lambda,1])$ equal to zero near
$\lambda$ and equal to one near $1$. For $\eta$ in a sufficiently
small open ball $\mathcal O_K\subset X_K$, define
\begin{equation}\label{eq:linear-pullback}
 T_\eta(r,\theta)=(s_\eta(r,\theta),\theta),\qquad
 s_\eta=r+\chi(r)\eta(\theta).
\end{equation}
The conditions $\partial_r s_\eta>0$ and $s_\eta>0$ hold when
$\|\eta\|_{C^{2,\alpha}}$ is small. Thus $T_\eta$ maps
$A_\lambda$ diffeomorphically onto
$\Omega_{1,\eta}\setminus\overline{B_\lambda}$ and leaves an
inner collar fixed. To give the transformed equation explicitly, set
\[
 a_\eta=1+\chi'\eta,\qquad b_\eta=\chi\eta_\theta,\qquad
 \mathcal D_\eta=a_\eta^{-1}\partial_r,\qquad
 \mathcal E_\eta=\partial_\theta-b_\eta a_\eta^{-1}\partial_r.
\]
The pullback of the Euclidean Laplacian is
\begin{equation}\label{eq:linear-pullback-op}
 \Delta_\eta v
 =\mathcal D_\eta^2v+s_\eta^{-1}\mathcal D_\eta v
     +s_\eta^{-2}\mathcal E_\eta^2v.
\end{equation}
Here $\mathcal D_\eta^2$ and $\mathcal E_\eta^2$ include the
derivatives of their coefficients. Formula
\eqref{eq:linear-pullback-op} shows that
$\eta\mapsto\Delta_\eta$ is $C^\infty$ with values in
$\mathcal L(C^{2,\alpha},C^{0,\alpha})$. These operators are
uniformly elliptic for small $\eta$ and preserve $D_K$ invariance.

\subsection{The radial branch and its Dirichlet nonresonance}
If $\mu\notin\Sigma_{\mathrm{rad}}(\lambda)$, the radial problem
\[
 U_\mu''+r^{-1}U_\mu'+\mu U_\mu=0,\qquad
 U_\mu(\lambda)=1,\quad U_\mu(1)=0
\]
has a unique solution. For $\mu>0$, write $t=\sqrt\mu$ and
\[
 \mathcal D_0(t)=J_0(t\lambda)Y_0(t)-Y_0(t\lambda)J_0(t).
\]
Then
\begin{equation}\label{eq:linear-U}
 U_\mu(r)=
 \frac{J_0(tr)Y_0(t)-Y_0(tr)J_0(t)}{\mathcal D_0(t)}.
\end{equation}
The denominator vanishes exactly at a radial Dirichlet eigenvalue.
The Bessel Wronskian identity gives
\begin{equation}\label{eq:linear-d}
 d_\mu:=U_\mu'(1)=-\frac{2}{\pi\mathcal D_0(\sqrt\mu)}\ne0.
\end{equation}
The nonvanishing also follows directly from uniqueness of the radial
Cauchy problem: $U_\mu(1)=U_\mu'(1)=0$ would imply
$U_\mu\equiv0$. We therefore have the radial solutions
$(\mu,0,Q_\mu)$ of \eqref{eq:linear}, where $Q_\mu=d_\mu^2>0$.
No positivity hypothesis on $U_\mu$ is needed for this problem.

We next verify that the hypothesis in Theorem~\ref{thm:linear}
indeed makes the full Dirichlet state equation invertible in the
chosen symmetry class. For $j\ge1$, let $\mu_{j,\lambda}^*$ be
the first eigenvalue of \eqref{eq:linear-Robin} with $K$ replaced
by $j$. Its Rayleigh characterization is
\begin{equation}\label{eq:linear-Rayleigh}
 \mu_{j,\lambda}^*
 =\min_{\substack{h\in H^1(\lambda,1)\setminus\{0\}\\h(\lambda)=0}}
 \frac{\displaystyle
   \int_\lambda^1\left(r|h'|^2+\frac{j^2}{r}|h|^2\right)\dd r
            +|h(1)|^2}
      {\displaystyle\int_\lambda^1r|h|^2\dd r}.
\end{equation}
In particular, $\mu_{j,\lambda}^*>0$. The first eigenfunction is
simple and may be chosen strictly positive on $(\lambda,1]$.
Comparing the quotients yields
\begin{equation}\label{eq:linear-mode-order}
 \mu_{j,\lambda}^*\ge\mu_{\ell,\lambda}^*+j^2-\ell^2
       \quad(j>\ell).
\end{equation}
Now we aim to establish radial nonresonance as follows.

\begin{lemma}\label{lem:linear-nonresonance}
At $\mu_*:=\mu_{K,\lambda}^*$, the map
\[
 -\Delta-\mu_*:Z_K^{\mathrm{lin}}\longrightarrow E_K^{\mathrm{lin}}
\]
is an isomorphism if and only if
$\mu_*\notin\Sigma_{\mathrm{rad}}(\lambda)$.
For $K=1$, this condition holds for every $\lambda\in(0,1)$.
For every fixed $K$, the exceptional set
\[
 E_K=\{\lambda\in(0,1):
        \mu_{K,\lambda}^*\in\Sigma_{\mathrm{rad}}(\lambda)\}
\]
is discrete in $(0,1)$ and is disjoint from a neighborhood of $1$.
\end{lemma}

\begin{proof}
We first compare the mixed and Dirichlet eigenvalues in each
nonzero angular mode. Let $\Lambda_j^D$ be the first Dirichlet
eigenvalue in mode $j\ge1$. Its Rayleigh characterization is
\[
 \Lambda_j^D
 =\min_{h\in H_0^1(\lambda,1)\setminus\{0\}}
 \frac{\displaystyle
   \int_\lambda^1\left(r|h'|^2+\frac{j^2}{r}|h|^2\right)\dd r}
      {\displaystyle\int_\lambda^1r|h|^2\dd r}.
\]
Every function in this form domain is admissible in
\eqref{eq:linear-Rayleigh}, and the boundary term $|h(1)|^2$
vanishes. Hence $\mu_{j,\lambda}^*\le\Lambda_j^D$.
To prove strictness, suppose that equality holds and choose a
nonzero Dirichlet minimizer $h$. It then also minimizes the mixed
quotient. Its Euler--Lagrange equation gives the same differential
equation and the natural outer condition $h'(1)+h(1)=0$.
Since $h(1)=0$, we also have $h'(1)=0$. Uniqueness of the
Cauchy problem at $r=1$ would give $h\equiv0$, a contradiction.
Thus $\Lambda_j^D>\mu_{j,\lambda}^*$, and
\eqref{eq:linear-mode-order} gives, for every permitted mode $j=nK$,
\[
 \Lambda_j^D>\mu_{j,\lambda}^*\ge\mu_{K,\lambda}^*=\mu_*.
\]

Let $w\in Z_K^{\mathrm{lin}}$ satisfy
$(-\Delta-\mu_*)w=0$. We denote its angular mean by $w_0(r)$
and its coefficient of $\cos(j\theta)$ by $w_j(r)$, $j=nK$.
The $D_K$ symmetries exclude all other nonzero Fourier modes.
Taking the corresponding coefficient in the equation gives
\[
 -w_j''-r^{-1}w_j'+j^2r^{-2}w_j=\mu_*w_j,\qquad
 w_j(\lambda)=w_j(1)=0.
\]
Multiplying by $rw_j$ and integrating by parts, we obtain
\[
 \mu_*\int_\lambda^1r|w_j|^2\dd r
 =\int_\lambda^1\left(r|w_j'|^2+\frac{j^2}{r}|w_j|^2\right)\dd r
 \ge\Lambda_j^D\int_\lambda^1r|w_j|^2\dd r.
\]
The strict inequality $\Lambda_j^D>\mu_*$ forces $w_j=0$.
The angular mean satisfies
\[
 -w_0''-r^{-1}w_0'=\mu_*w_0,\qquad
 w_0(\lambda)=w_0(1)=0.
\]
This radial problem has only the zero solution precisely when
$\mu_*$ is outside the radial Dirichlet spectrum. Conversely, every radial Dirichlet eigenfunction at
$\mu_*$ belongs to $Z_K^{\mathrm{lin}}$ and gives a nonzero
kernel element. Uniqueness of Fourier coefficients therefore shows
that the kernel is trivial if and only if
$\mu_*\notin\Sigma_{\mathrm{rad}}(\lambda)$.

We next pass from this kernel statement to invertibility. The
Dirichlet Laplacian is an isomorphism
$Z_K^{\mathrm{lin}}\to E_K^{\mathrm{lin}}$ by Schauder theory
\cite{GT}. Indeed, the usual Dirichlet solution for a
$D_K$-invariant right-hand side is $D_K$-invariant by uniqueness.
Let $J_D:Z_K^{\mathrm{lin}}\to E_K^{\mathrm{lin}}$ be the
natural inclusion. The compact embedding
$C^{2,\alpha}(\overline A_\lambda)\hookrightarrow
C^{0,\alpha}(\overline A_\lambda)$ makes $J_D$ compact.
We have
\[
 (-\Delta)^{-1}(-\Delta-\mu_*)
 =\Id_{Z_K^{\mathrm{lin}}}-\mu_* (-\Delta)^{-1}J_D.
\]
The right-hand side is the identity plus a compact operator, so it
is Fredholm of index zero. Thus $-\Delta-\mu_*$ has index zero
as well. A trivial kernel then gives surjectivity and, by the
bounded inverse theorem, an isomorphism between the stated spaces.
This proves the first assertion.

For $K=1$, we can exclude radial resonance directly. Choose a
positive first mixed eigenfunction $h$, set
$I=\int_\lambda^1h(t)\dd t>0$, and define
\[
 W(r)=I^{-1}\int_r^1h(t)\dd t.
\]
The equation for $h$ can be written as
\[
 (h'+h/r)'=-\mu_{1,\lambda}^*h.
\]
Integrating from $r$ to $1$ and using $h'(1)+h(1)=0$, we get
\[
 h'(r)+\frac{h(r)}r
    =\mu_{1,\lambda}^*\int_r^1h(t)\dd t.
\]
Since $W'=-h/I$ and $W''=-h'/I$, it follows that
\[
 -W''-r^{-1}W'=\frac{h'+h/r}{I}
                  =\mu_{1,\lambda}^*W,\qquad
 W'(\lambda)=0,\quad W(1)=0.
\]
Moreover, $W(\lambda)=1$ and $W>0$ on $[\lambda,1)$.
The positivity of $W$ identifies $\mu_{1,\lambda}^*$ as the
first radial Neumann--Dirichlet eigenvalue. To see this directly,
let $v$ be a positive first eigenfunction of that problem, with
eigenvalue $\nu_{ND}$. Green's identity and the common boundary
conditions give
\[
 (\mu_{1,\lambda}^*-\nu_{ND})
       \int_\lambda^1rWv\dd r=0.
\]
The integral is positive, so $\mu_{1,\lambda}^*=\nu_{ND}$.

The corresponding Rayleigh quotient is
\[
 \nu_{ND}
 =\min_{\substack{v\in H^1(\lambda,1)\setminus\{0\}\\v(1)=0}}
       \frac{\displaystyle\int_\lambda^1r|v'|^2\dd r}
            {\displaystyle\int_\lambda^1r|v|^2\dd r}.
\]
Adding the condition $v(\lambda)=0$ gives the first radial
Dirichlet eigenvalue, which we denote by $\nu_D$.
Hence $\nu_{ND}\le\nu_D$. If equality held, a Dirichlet
minimizer would also minimize the Neumann--Dirichlet quotient and
would satisfy its natural inner condition $v'(\lambda)=0$.
Together with $v(\lambda)=0$, this would force $v=0$ by
Cauchy uniqueness. Thus
$\mu_{1,\lambda}^*=\nu_{ND}<\nu_D$. Since every radial
Dirichlet eigenvalue is at least $\nu_D$, we have $E_1=\varnothing$.

We now verify the absence of exceptional radii near $1$ for every
fixed $K$. Put $\delta=1-\lambda$ and use
$h(r)=\sin(\pi(r-\lambda)/(2\delta))$ in
\eqref{eq:linear-Rayleigh}. The elementary integrals are
\[
 \int_\lambda^1h^2\dd r=\frac\delta2,\qquad
 \int_\lambda^1|h'|^2\dd r=\frac{\pi^2}{8\delta},\qquad
 h(1)=1.
\]
Since $\lambda\le r\le1$, we have
\[
 \int_\lambda^1rh^2\dd r\ge\frac{\lambda\delta}{2},\qquad
 \int_\lambda^1r|h'|^2\dd r\le\frac{\pi^2}{8\delta},\qquad
 \int_\lambda^1\frac{K^2}{r}h^2\dd r
       \le\frac{K^2\delta}{2\lambda}.
\]
Substituting these bounds into the quotient gives
\begin{equation}\label{eq:thin-linear-upper}
 \mu_{K,\lambda}^*
 \le\frac{\pi^2}{4\lambda\delta^2}
       +\frac{K^2}{\lambda^2}+\frac{2}{\lambda\delta}.
\end{equation}
For a function $v$ with zero Dirichlet trace on both circles, the
one-dimensional Poincar\'e inequality on each radial ray gives
\[
 \int_\lambda^1r|\partial_r v|^2\dd r
 \ge\lambda\int_\lambda^1|\partial_r v|^2\dd r
 \ge\frac{\lambda\pi^2}{\delta^2}\int_\lambda^1|v|^2\dd r
 \ge\frac{\lambda\pi^2}{\delta^2}\int_\lambda^1r|v|^2\dd r.
\]
Integrating in $\theta$ and adding the nonnegative angular part
of the energy gives, for the first full Dirichlet eigenvalue,
\begin{equation}\label{eq:thin-linear-lower}
 \Lambda_D(\lambda)\ge\frac{\lambda\pi^2}{\delta^2}.
\end{equation}
For fixed $K$, the upper bound for
$\delta^2\mu_{K,\lambda}^*$ tends to $\pi^2/4$ as
$\lambda\to1$, while the lower bound for
$\delta^2\Lambda_D(\lambda)$ tends to $\pi^2$.
Thus $\mu_{K,\lambda}^*<\Lambda_D(\lambda)$ for all
$\lambda$ sufficiently close to $1$. Every radial Dirichlet
eigenvalue is at least $\Lambda_D(\lambda)$, so these radii
lie outside $E_K$.

Finally, we prove that $E_K$ is discrete in $(0,1)$.
Let $h=h(r;\lambda,\mu)$ solve
\[
 -h''-r^{-1}h'+K^2r^{-2}h=\mu h,\qquad
 h(\lambda)=0,\quad h'(\lambda)=1,
\]
and set $B(\lambda,\mu)=h'(1)+h(1)$.
After mapping $[\lambda,1]$ to a fixed interval, the coefficients
and initial data are real analytic in $(\lambda,\mu)$.
Analytic dependence for the Cauchy problem therefore makes $B$
real analytic. The value $\mu_{K,\lambda}^*$ is its principal
mixed root. With $\lambda$ fixed, write $\dot h=\partial_\mu h$.
Differentiating the equation and initial data gives
\[
 -(r\dot h')'+\frac{K^2}{r}\dot h-\mu r\dot h=rh,\qquad
 \dot h(\lambda)=\dot h'(\lambda)=0.
\]
Multiplying this equation by $h$ and subtracting $\dot h$ times
the equation for $h$, we obtain, at a mixed root,
\[
 \int_\lambda^1rh^2\dd r
 =\bigl[r(\dot h h'-h\dot h')\bigr]_\lambda^1
 =-h(1)\bigl(\dot h'(1)+\dot h(1)\bigr).
\]
Here the inner term vanishes by the initial data, and the last
identity uses $h'(1)=-h(1)$. Consequently,
\[
 h(1)\partial_\mu B(\lambda,\mu)
       =-\int_\lambda^1rh(r)^2\dd r\ne0.
\]
In particular, $\partial_\mu B\ne0$ at the principal root.
The analytic implicit function theorem gives a local analytic
branch of mixed eigenvalues. By simplicity and continuity of the
first eigenvalue, these local branches coincide with
$\mu_{K,\lambda}^*$ and agree on overlaps. Hence
$\lambda\mapsto\mu_{K,\lambda}^*$ is real analytic on $(0,1)$.

Since $\mu_{K,\lambda}^*>0$, the function
\[
 F_K(\lambda)=\mathcal D_0(\sqrt{\mu_{K,\lambda}^*})
\]
is real analytic as well. Here
\[
 \mathcal D_0(t)=J_0(t\lambda)Y_0(t)-Y_0(t\lambda)J_0(t),
\]
so its dependence on both $\lambda$ and $t>0$ is real analytic.
The Bessel determinant vanishes exactly at radial Dirichlet
eigenvalues, so the zero set of $F_K$ is $E_K$.
The thin-annulus comparison shows that $F_K$ is nonzero near $1$,
and hence it is not identically zero. The zeros of a nonzero real
analytic function on an interval are isolated. This proves the
last assertion.
\end{proof}

\subsection{The state equation and the Bernoulli operator}
Assume henceforth the hypothesis of Theorem~\ref{thm:linear}, and
write $\mu_*=\mu_{K,\lambda}^*$. We now construct the actual
operator equation for the moving-boundary problem. Set
$b(r)=\log r/\log\lambda$, so that $b(\lambda)=1$ and $b(1)=0$.
For $w\in Z_K^{\mathrm{lin}}$, define
\begin{equation}\label{eq:linear-state-map}
 \mathcal F_{\mathrm{lin}}(\mu,\eta,w)
       =-\Delta_\eta(b+w)-\mu(b+w)
       \in E_K^{\mathrm{lin}}.
\end{equation}
The parameters range over an open neighborhood of
$(\mu_*,0,U_{\mu_*}-b)$ in
$\R\times X_K\times Z_K^{\mathrm{lin}}$. By
\eqref{eq:linear-pullback-op}, this map is $C^\infty$, and
\[
 D_w\mathcal F_{\mathrm{lin}}(\mu_*,0,U_{\mu_*}-b)
             =-\Delta-\mu_*.
\]

\begin{proposition}
\label{prop:linear-state}
There are an interval $I\subset(0,\infty)$ containing $\mu_*$ and an open ball
$\mathcal O_K\subset X_K$ containing zero such that the Dirichlet
problem
\begin{equation}\label{eq:linear-moving-state}
 \begin{cases}
  \Delta u+\mu u=0&\text{in }\Omega_{1,\eta}\setminus\overline{B_\lambda},\\
  u=1&\text{on }\partial B_\lambda,\\
  u=0&\text{on }\partial\Omega_{1,\eta}
 \end{cases}
\end{equation}
has a unique $D_K$-invariant solution $u_{\mu,\eta}$ for
$(\mu,\eta)\in I\times\mathcal O_K$. Its pullback
$v_{\mu,\eta}=u_{\mu,\eta}\circ T_\eta$ belongs to
$b+Z_K^{\mathrm{lin}}$ and depends $C^\infty$ on $(\mu,\eta)$.
At $\eta=0$, it equals $U_\mu$.
\end{proposition}

\begin{proof}
Lemma~\ref{lem:linear-nonresonance} makes $D_w\mathcal F_{\mathrm{lin}}$
an isomorphism at the base point. The implicit function theorem
therefore gives $w=w(\mu,\eta)$ and the smooth pullback
$v_{\mu,\eta}=b+w(\mu,\eta)$. Moreover,
$-\Delta_\eta-\mu$ remains an isomorphism on the same zero-trace
spaces after shrinking the neighborhoods, by continuity in operator
norm. Thus uniqueness holds for the entire linear Dirichlet problem
in the symmetry class, not merely among nearby states. At $\eta=0$,
the radial formula \eqref{eq:linear-U} provides that solution.
\end{proof}

Let $\Pi g=g-\av g$ and write the outward unit normal to the graph
$r=1+\eta(\theta)$ as
\[
 \nu_\eta
  =\frac{(1+\eta)e_r-\eta_\theta e_\theta}
         {\sqrt{(1+\eta)^2+\eta_\theta^2}}.
\]
We define
\begin{equation}\label{eq:linear-N-and-G}
 \begin{split}
 N_{\mathrm{lin}}(\mu,\eta)(\theta)
   &=\partial_{\nu_\eta}u_{\mu,\eta}
                   ((1+\eta(\theta))e_r),\\
 G_{\mathrm{lin}}(\mu,\eta)
   &=\Pi\bigl(N_{\mathrm{lin}}(\mu,\eta)^2\bigr).
 \end{split}
\end{equation}
Since $\chi=1$ near $1$, the constant outer Dirichlet trace gives
the explicit fixed-boundary formula
\begin{equation}\label{eq:linear-normal-trace}
 N_{\mathrm{lin}}(\mu,\eta)
   =\left(1+\frac{\eta_\theta^2}{(1+\eta)^2}\right)^{1/2}
                      \partial_r v_{\mu,\eta}(1,\cdot).
\end{equation}
Consequently,
\[
 N_{\mathrm{lin}}:I\times\mathcal O_K
           \longrightarrow C^{1,\alpha}_{D_K}(\Sone),\qquad
 G_{\mathrm{lin}}:I\times\mathcal O_K\longrightarrow Y_K
\]
are $C^\infty$ maps. At a radial state,
$N_{\mathrm{lin}}(\mu,0)=d_\mu$ and
$G_{\mathrm{lin}}(\mu,0)=0$. After shrinking the neighborhoods,
$N_{\mathrm{lin}}$ is nowhere zero because $d_{\mu_*}\ne0$.
Since the tangential derivative of $u_{\mu,\eta}$ vanishes on the
outer boundary, the full overdetermined problem is equivalent to
\begin{equation}\label{eq:recover-Q}
 G_{\mathrm{lin}}(\mu,\eta)=0,\qquad
 Q=\av{N_{\mathrm{lin}}(\mu,\eta)^2}>0.
\end{equation}
The mean is taken in the angular coordinate; its role is only to
remove constants, so no arclength weight is needed.

\subsection{Linearization, kernel, and range}\label{sec:boundary}
Fix $\mu\in I$ and $\xi\in X_K$. The material derivative of
$v_{\mu,\varepsilon\xi}$ has zero Dirichlet trace. Its Eulerian
counterpart is
\[
 u_\xi'=\left.\partial_\varepsilon
                 v_{\mu,\varepsilon\xi}\right|_{\varepsilon=0}
                       -\chi(r)U_\mu'(r)\xi(\theta).
\]
It follows from radial ODE that $U_\mu''(1)=-U_\mu'(1)=-d_\mu$.
Then differentiating \eqref{eq:linear-moving-state}, as in the standard
shape-derivative calculation \cite[Lemma 3.2]{GWYZ}, we have
\begin{equation}\label{eq:linear-shape-state}
 \begin{cases}
  (-\Delta-\mu)u_\xi'=0&\text{in }A_\lambda,\\
  u_\xi'=0&\text{on }r=\lambda,\\
  u_\xi'=-U_\mu'(1)\xi=-d_\mu\xi&\text{on }r=1.
 \end{cases}
\end{equation}
 Either differentiating
\eqref{eq:linear-normal-trace} or differentiating the physical normal
trace therefore yields
\begin{equation}\label{eq:linear-shape-normal}
 D_\eta N_{\mathrm{lin}}(\mu,0)\xi=\partial_r u_\xi'(1,\cdot)+U_\mu''(1)\xi
       =\partial_r u_\xi'(1,\cdot)-d_\mu\xi.
\end{equation}
The normal-vector variation is tangential and has zero scalar product
with the radial gradient of $U_\mu$.

For every permitted frequency $j=nK$, let $H_{j,\mu}$ solve
\[
 -H''-r^{-1}H'+j^2r^{-2}H-\mu H=0,\qquad
 H(\lambda)=0,\quad H(1)=1.
\]
The Dirichlet invertibility established above makes this solution
well defined. In Bessel notation it is
\begin{equation}\label{eq:linear-H}
 H_{j,\mu}(r)=
 \frac{J_j(tr)Y_j(t\lambda)-Y_j(tr)J_j(t\lambda)}
      {J_j(t)Y_j(t\lambda)-Y_j(t)J_j(t\lambda)},\qquad t=\sqrt\mu.
\end{equation}
For $\xi=\cos(j\theta)$, equation
\eqref{eq:linear-shape-state} gives
$u_\xi'=-d_\mu H_{j,\mu}(r)\cos(j\theta)$. Thus, with
\begin{equation}\label{eq:linear-symbol}
 T_{j,\lambda}(\mu)=H_{j,\mu}'(1)+1,\qquad
 m_j^{\mathrm{lin}}(\mu)=-2d_\mu^2T_{j,\lambda}(\mu),
\end{equation}
we have
\begin{equation}\label{eq:linear-multiplier}
 L_\mu[\cos(j\theta)]
 :=D_\eta G_{\mathrm{lin}}(\mu,0)[\cos(j\theta)]
 =m_j^{\mathrm{lin}}(\mu)\cos(j\theta).
\end{equation}
The symbol $T_{j,\lambda}$ vanishes precisely when
\[
 -h''-r^{-1}h'+j^2r^{-2}h=\mu h,\qquad
 h(\lambda)=0,\quad h'(1)+h(1)=0
\]
has a nonzero solution. In particular, at $\mu=\mu_*$ the
$K$-mode is critical. By \eqref{eq:linear-mode-order}, all the
mixed eigenvalues in the modes $j=nK$, $n\ge2$, lie strictly
above $\mu_*$. Hence $T_{nK,\lambda}(\mu_*)\ne0$ for $n\ge2$.
The following lemma establishes the Fredholm property and characterizes the range for the linear law.

\begin{lemma}
\label{lem:Fredholm}
For $\mu\in I$, the operator $L_\mu:X_K\to Y_K$ is Fredholm
of index zero. At $\mu=\mu_*$,
\begin{equation}\label{eq:linear-kernel-range}
 \begin{split}
 \Ker L_{\mu_*}&=\Span\{\cos(K\theta)\},\\
 \Ran L_{\mu_*}&=\left\{g\in Y_K:
          \int_0^{2\pi}g(\theta)\cos(K\theta)\dd\theta=0\right\}.
 \end{split}
\end{equation}
\end{lemma}

\begin{proof}
For $\xi\in X_K$, let $P_\mu\xi$ be the unique $D_K$-invariant
solution of
\[
 (-\Delta-\mu)P_\mu\xi=0\quad\text{in }A_\lambda,\qquad
 P_\mu\xi|_{r=\lambda}=0,\quad P_\mu\xi|_{r=1}=\xi.
\]
The Dirichlet invertibility on $I$ and the Schauder estimates
give a bounded extension $P_\mu:X_K\to
C^{2,\alpha}_{D_K}(\overline A_\lambda)$. Its angular mean
solves the homogeneous radial Dirichlet problem, so it vanishes
by uniqueness. Consequently,
\[
 \Lambda_\mu\xi=\partial_r(P_\mu\xi)(1,\cdot)
\]
defines a bounded map $X_K\to Y_K$. The shape-derivative formulas
give
\[
 L_\mu=-2d_\mu^2 A_\mu,\qquad A_\mu=\Lambda_\mu+J,
\]
where $J:X_K\hookrightarrow Y_K$ is the natural inclusion.
This is a Dirichlet-to-Neumann operator plus a zero-order term.
By the standard elliptic theory for Dirichlet-to-Neumann
operators, cf. the proof of
\cite[Proposition~3.2]{Sicbaldi}, $L_\mu$ is Fredholm.

We use the corresponding Dirichlet--Robin problem to recall its
Fredholm alternative in the precise spaces needed here:
\[
 \begin{cases}
  (-\Delta-\mu)v=0&\text{in }A_\lambda,\\
  v=0&\text{on }r=\lambda,\\
  \partial_r v+v=g&\text{on }r=1.
 \end{cases}
\]
We restrict to $D_K$-invariant interior functions with zero angular
mean at each radius. On this $H^1$ space, the weak form is symmetric and becomes
coercive after adding a sufficiently large $L^2$ term. The compact
embedding into $L^2$ gives the self-adjoint Fredholm alternative:
the homogeneous solution space is finite dimensional, and the
problem is solvable precisely when
\[
 \int_0^{2\pi}g(\theta)v_0(1,\theta)\dd\theta=0
\]
for every homogeneous solution $v_0$. For $g\in Y_K$, the mixed
boundary Schauder estimates give $v\in C^{2,\alpha}$; the two
boundary components are smooth and disjoint \cite{GT}.
Dirichlet uniqueness identifies these homogeneous outer traces
with $\Ker A_\mu$, and the outer trace of a solution satisfies
$A_\mu[v(1,\cdot)]=g$. Hence
\[
 \Ran A_\mu=\left\{g\in Y_K:
     \int_0^{2\pi}g\phi\dd\theta=0
                  \text{ for every }\phi\in\Ker A_\mu\right\}.
\]
The finite-dimensional kernel lies in $X_K\subset Y_K$, so these
orthogonality conditions are independent. Thus the range is closed
and has codimension $\dim\Ker A_\mu$. This proves that $A_\mu$,
and therefore $L_\mu$ since $d_\mu\ne0$, is Fredholm of index zero.
The same argument applies to any real radial $C^{0,\alpha}$
potential and a constant Robin coefficient whenever the Dirichlet
problem is invertible in the chosen symmetry class.

At $\mu_*$, taking Fourier coefficients in the state equation and
using \eqref{eq:linear-multiplier} shows that the coefficient in
mode $j=nK$ of $L_{\mu_*}\xi$ equals
$m_j^{\mathrm{lin}}(\mu_*)$ times that of $\xi$.
The multiplier vanishes for $n=1$ and is nonzero for every $n\ge2$,
so uniqueness of Fourier coefficients gives
$\Ker L_{\mu_*}=\Span\{\cos(K\theta)\}$.
The continuous functional
\[
 \ell(g)=\int_0^{2\pi}g(\theta)\cos(K\theta)\dd\theta
\]
vanishes on every range element and satisfies
$\ell(\cos(K\theta))=\pi\ne0$. Thus the range is contained
in the closed codimension-one subspace $\Ker\ell$.
Since $L_{\mu_*}$ is Fredholm of index zero and has a
one-dimensional kernel, its range is closed and has codimension one
as well. The inclusion is therefore an equality, which proves
\eqref{eq:linear-kernel-range}. This argument does not require
Fourier partial sums to converge in the $C^{2,\alpha}$ norm.
\end{proof}

\subsection{Transversality and proof of Theorem~\ref{thm:linear}}
We verify the remaining Crandall--Rabinowitz condition explicitly.
A dot will denote differentiation with respect to $\mu$. Since
$H_{K,\mu}(\lambda)=0$ and $H_{K,\mu}(1)=1$ for all $\mu$,
\[
 \left(-\frac{d^2}{dr^2}-\frac1r\frac d{dr}
                +\frac{K^2}{r^2}-\mu\right)\dot H_{K,\mu}
       =H_{K,\mu},\qquad
 \dot H_{K,\mu}(\lambda)=\dot H_{K,\mu}(1)=0.
\]
Multiplying the equation for $\dot{H}_{K,\mu}$ by $rH$ and the equation for $H_{K,\mu}$ by $r\dot{H}$ and subtracting the latter from the former, we can obtain that
\[
\int_{\lambda}^1rH^2_{K,\mu}dr=-\dot  H'_{K,\mu}(1),
\]
which gives
\begin{equation}\label{eq:linear-crossing}
\frac{d T_{K,\lambda}(\mu)}{d \mu}
        =-\int_\lambda^1rH_{K,\mu}(r)^2\dd r<0.
\end{equation}
Because $T_{K,\lambda}(\mu_*)=0$, differentiating
\eqref{eq:linear-symbol} at the critical value yields
\begin{equation}\label{eq:linear-CR-crossing}
\frac{d m_K^{\mathrm{lin}}}{d\mu}(\mu_*)
      =2d_{\mu_*}^2\int_\lambda^1rH_{K,\mu_*}(r)^2\dd r>0.
\end{equation}
In particular,
\[
 D_\mu D_\eta G_{\mathrm{lin}}(\mu_*,0)\cos(K\theta)
    \notin\Ran D_\eta G_{\mathrm{lin}}(\mu_*,0).
\]
We have now checked all hypotheses on the concrete map
$G_{\mathrm{lin}}:I\times\mathcal O_K\to Y_K$: it is smooth,
$G_{\mathrm{lin}}(\mu,0)=0$, its derivative in $\eta$ has index
zero and a one-dimensional kernel, and the crossing is transverse.
The Crandall--Rabinowitz theorem \cite{CR} gives a smooth local
curve of nontrivial solutions. Choose the parameter by
\[
 s=\frac1\pi\int_0^{2\pi}\eta(s,\theta)\cos(K\theta)\dd\theta.
\]
With the complement of the kernel defined by vanishing $K$-th
coefficient, the curve has
\[
 \mu(0)=\mu_*,\qquad
 \eta(s)=s\cos(K\theta)+O(s^2)\quad\text{in }X_K.
\]
Every nearby nontrivial solution belongs to this curve.

Let $\mathcal R\eta(\theta)=\eta(\theta+\pi/K)$. This
half-period rotation preserves $X_K$, $Y_K$, and the chosen
complement, and reverses the $K$-th Fourier coefficient. Rotation
invariance of the state equation and uniqueness of its Dirichlet
solution give
$G_{\mathrm{lin}}(\mu,\mathcal R\eta)
=\mathcal R G_{\mathrm{lin}}(\mu,\eta)$.
The equivariance argument \cite{GSS} and local uniqueness of the
bifurcation curve therefore yield
\[
 \mu(-s)=\mu(s),\qquad \eta(-s)=\mathcal R\eta(s).
\]
Hence $\mu'(0)=0$, and smoothness gives
$\mu(s)=\mu_*+O(s^2)$. Formula \eqref{eq:recover-Q} defines
the corresponding smooth positive $Q(s)$. This completes the proof
of Theorem~\ref{thm:linear}.

\section{The power-law overdetermined problem}\label{sec:nonlinear-proof}

We now prove Theorem~\ref{thm:power} for the specific equation
$\Delta u=u-u^p$. Throughout this section, $\lambda>0$ and $p>1$
are fixed, $a=2/(p-1)$, and $V_p=V$ denotes the profile in
\eqref{eq:pulse-intro}. We use the mean outer radius $R$ as the
parameter. We first specify the spaces for this problem and construct
its large radial branch. We then solve its Dirichlet equation locally,
form its Bernoulli equation, and verify the kernel, range, and
transversality conditions for the Crandall--Rabinowitz theorem.

\subsection{The working spaces and the power map}
Fix temporarily a reference radius $R_*>\lambda$ and put
$A_*=B_{R_*}\setminus\overline{B_\lambda}$. For an integer $K\ge2$,
let $D_K$ be generated by $\theta\mapsto-\theta$ and
$\theta\mapsto\theta+2\pi/K$. We retain the boundary spaces
$X_K$ and $Y_K$ from the introduction and use the ordinary H\"older spaces
\[
 Z_{p,K}=\{w\in C^{2,\alpha}_{D_K}(\overline A_*):
             w|_{r=\lambda}=w|_{r=R_*}=0\},\qquad
 W_{p,K}=C^{0,\alpha}_{D_K}(\overline A_*).
\]
Every $D_K$-invariant state with the prescribed Dirichlet values can be written as
$v=b_*+w$, where $w\in Z_{p,K}$ and
\[
 b_*(r)=\frac{R_*-r}{R_*-\lambda}.
\]
We work in a sufficiently small $C^{2,\alpha}$ neighborhood of a
positive radial state. The radial mode remains in the state and
equation spaces. The zero-mean condition applies only to the boundary
perturbation, since $R$ already describes the mean outer radius.

For nonintegral $p$, we need to check smoothness of the power map
at a state with zero outer trace. Put $\varrho(r)=R_*-r$.
For every $v\in C^{2,\alpha}(\overline A_*)$ vanishing at $R_*$,
we have, with $t=R_*-r$,
\begin{equation}\label{eq:division}
 \frac{v(t,\theta)}t=\int_0^1\partial_t v(\tau t,\theta)\dd\tau.
\end{equation}
Thus the quotient extends across $t=0$ and satisfies
\[
 \|v/\varrho\|_{C^{1,\alpha}(\overline A_*)}
       \le C\|v\|_{C^{2,\alpha}(\overline A_*)}.
\]
This estimate allows us to use ordinary small H\"older neighborhoods
without introducing a separate space for the quotient.
We first describe the power map in the Dirichlet space.

\begin{lemma}\label{lem:power-map}
Let $U\in C^{2,\alpha}_{D_K}(\overline A_*)$ be positive in $A_*$,
with $U|_{r=\lambda}=1$, $U|_{r=R_*}=0$, and
$\partial_rU(R_*,\theta)<0$. For a sufficiently small
$\epsilon_U>0$, the map
\[
 w\longmapsto(U+w)^p:
 \{w\in Z_{p,K}:\|w\|_{C^{2,\alpha}}<\epsilon_U\}
        \longrightarrow W_{p,K}
\]
is $C^\infty$. Every state $U+w$ in this neighborhood is positive
in $A_*$ and has negative outer normal derivative.
The same assertion holds without imposing $D_K$ invariance.
\end{lemma}

\begin{proof}
The quotient $U/\varrho$ is positive on the closure: its boundary
value at $R_*$ is $-\partial_rU(R_*,\theta)>0$, its inner value
is $1/(R_*-\lambda)$, and $U>0$ in the interior. By compactness,
\[
 c_U:=\min_{\overline A_*}U/\varrho>0.
\]
The division estimate gives
$\|w/\varrho\|_{C^{1,\alpha}}\le C\|w\|_{C^{2,\alpha}}$.
Choosing $C\epsilon_U<c_U/2$, we obtain
\[
 z:=\frac{U+w}{\varrho}\ge c_U/2>0.
\]
This proves positivity in the annulus and, at the outer boundary,
$\partial_r(U+w)=-z<0$.

The map $w\mapsto z$ is affine and continuous into the Banach
algebra $C^{1,\alpha}$. The scalar function $z\mapsto z^p$
is smooth on the set bounded away from zero, so it gives a
$C^\infty$ map on this neighborhood. Since
$\varrho^p\in C^{0,\alpha}$ for $p>1>\alpha$, multiplication
by $\varrho^p$ is bounded into $C^{0,\alpha}$. Hence
$(U+w)^p=\varrho^p z^p$ is $C^\infty$ with the asserted target.
Its derivatives in zero-trace directions $w_i\in Z_{p,K}$ are
\[
 D^m(v^p)[w_1,\ldots,w_m]
 =p(p-1)\cdots(p-m+1)\varrho^p z^{p-m}
       \prod_{i=1}^m\frac{w_i}{\varrho},\qquad v=U+w.
\]
Each quotient belongs to $C^{1,\alpha}$, and $z$ stays uniformly
positive. The formula therefore defines bounded multilinear maps
into $C^{0,\alpha}$ at every order, with continuous dependence on
$v$. The boundary vanishing of the directions supplies the factors
of $\varrho$ that remove the apparent singularity of $v^{p-m}$.
The symmetry is preserved by taking powers. The proof without
symmetry is identical.
\end{proof}

\subsection{The trivial branch: large positive radial solutions}\label{sec:radial}

Fix $\lambda>0$ and write $R=\lambda+\delta$, where
$\delta>0$ is sufficiently small. We seek positive radial
solutions of
\begin{equation}\label{eq:radial-original}
 U''+\frac1rU'=U-U^p,\qquad U(\lambda)=1,\quad U(\lambda+\delta)=0.
\end{equation}
We shall construct a family of size $\delta^{-a}$, where
$a=2/(p-1)$. This large family is distinct from any small ordered
solutions of the Dirichlet problem. In this construction, we use the
positive solution $V_p$ of \eqref{eq:pulse-intro} and the
Dirichlet nondegeneracy in Lemma~\ref{lem:pulse-D}. We prove these
one-dimensional facts in the next subsection directly from the
limiting equation, independently of the annular construction.

\subsubsection{Rescaling and the implicit-function theorem}
The radial second derivative has size $\delta^{-a-2}$ when a
function of size $\delta^{-a}$ varies across an interval of length
$\delta$, whereas its $p$th power has size $\delta^{-ap}$.
To retain both terms in the limit, we require $a+2=ap$, which
gives $a=2/(p-1)$. Accordingly, we set
\begin{equation}\label{eq:radial-scaling}
 t=\frac{r-\lambda}{\delta},\qquad
 U(r)=\delta^{-a}V_\delta(t).
\end{equation}
The derivatives and the power are
\[
 U'(r)=\delta^{-a-1}V_\delta'(t),\qquad
 U''(r)=\delta^{-a-2}V_\delta''(t),\qquad
 U(r)^p=\delta^{-ap}V_\delta(t)^p.
\]
We substitute these expressions into \eqref{eq:radial-original}
and multiply by $\delta^{a+2}$. Since $a+2-ap=0$, the
second derivative and the power retain coefficient $1$, the radial
first derivative has coefficient $\delta/(\lambda+\delta t)$,
and the linear term has coefficient $\delta^2$. The boundary
values become $V_\delta(0)=\delta^aU(\lambda)=\delta^a$
and $V_\delta(1)=0$. Thus the scaled problem is
\begin{equation}\label{eq:Vdelta}
 \begin{cases}
 V_\delta''+\dfrac{\delta}{\lambda+\delta t}V_\delta'
       -\delta^2V_\delta+V_\delta^p=0,&0<t<1,\\
 V_\delta(0)=\delta^a,\qquad V_\delta(1)=0.
 \end{cases}
\end{equation}
As the width tends to zero, the curvature coefficient and the linear
reaction term vanish. The rescaled inner trace also tends to zero,
although the original inner value remains $1$. These observations
give the limiting equation \eqref{eq:pulse-intro} and its two zero
endpoint values.

We cannot apply the smooth implicit function theorem directly at
$\delta=0$ with the trace $\delta^a$, because $a$ need not be an
integer and this trace need not even be differentiable at zero.
Instead, we introduce an independent real boundary parameter $b$
and write $V=v+b(1-t)$, with $v(0)=v(1)=0$. On
$C^{2,\alpha}_0([0,1])=\{v\in C^{2,\alpha}([0,1]):v(0)=v(1)=0\}$,
we define the residual with values in $C^{0,\alpha}([0,1])$ by
\begin{align}
 \mathcal E(\delta,b,v)
 ={}&v''+\frac{\delta}{\lambda+\delta t}\bigl(v'-b\bigr)
       -\delta^2\bigl(v+b(1-t)\bigr) \notag\\
    &+\bigl(v+b(1-t)\bigr)_+^p.
 \label{eq:limit-residual}
\end{align}
Here $s_+=\max\{s,0\}$. For $|\delta|<\lambda/2$, the
denominator is positive on $[0,1]$. The positive-part extension
therefore defines the residual on an open neighborhood of
$(0,0,V_p)$, including parameters of both signs. We shall use only
$b=\delta^a>0$ and prove that the resulting profiles are positive.
The following lemma constructs the radial family.

\begin{lemma}\label{lem:radial-family}
For every sufficiently small $\delta>0$, problem \eqref{eq:Vdelta}
admits a positive solution $V_\delta$ satisfying
\begin{equation}\label{eq:V-estimate}
 \|V_\delta-V_p\|_{C^{2,\alpha}([0,1])}
 \le C(\delta+\delta^a).
\end{equation}
The family is locally unique among scaled profiles close to $V_p$,
is $C^\infty$ for $\delta>0$, and satisfies
\begin{equation}\label{eq:V-derivative}
 \|\delta\,\partial_\delta V_\delta\|_{C^{2,\alpha}([0,1])}
 \le C(\delta+\delta^a).
\end{equation}
\end{lemma}

\begin{proof}
We first verify the regularity needed at $(0,0,V_p)$. Let
$F(v)=(v_+)^p$ and choose
$\alpha<\beta<\min\{1,p-1\}$. The scalar derivative
$F'(s)=p(s_+)^{p-1}$ is H\"older continuous with every such
exponent $\beta$ on bounded intervals. Hence, on a bounded set of
$C^1$ functions, the coefficients $p(v_+)^{p-1}$ are uniformly
bounded in $C^{0,\beta}$. If $v_n\to v$ in $C^{2,\alpha}$,
these coefficients converge uniformly by continuity of the scalar
power. The interpolation inequality
\[
 \|g\|_{C^{0,\alpha}}
 \le C\|g\|_{C^0}^{1-\alpha/\beta}
          \|g\|_{C^{0,\beta}}^{\alpha/\beta}
\]
then gives their convergence in $C^{0,\alpha}$. Multiplication
by such a coefficient is bounded from $C^{2,\alpha}$ into
$C^{0,\alpha}$, so the candidate derivative
$DF(v)w=p(v_+)^{p-1}w$ varies continuously in operator norm.
Finally, the pointwise integral formula
\[
 F(v+w)-F(v)-p(v_+)^{p-1}w
 =\int_0^1p\bigl((v+sw)_+^{p-1}-(v_+)^{p-1}\bigr)w\dd s
\]
shows that the remainder is $o(\|w\|_{C^{2,\alpha}})$ in
$C^{0,\alpha}$. Thus $F$ is $C^1$ between these spaces, and
the same is true of $\mathcal E$ near $(0,0,V_p)$.

Since $V_p$ solves \eqref{eq:pulse-intro}, we have
$\mathcal E(0,0,V_p)=0$. The derivative in the zero-trace
variable is
\[
 D_v\mathcal E(0,0,V_p)w
 =w''+pV_p^{p-1}w=-A_Dw,\qquad w(0)=w(1)=0.
\]
Lemma~\ref{lem:pulse-D} gives $\Ker A_D=\{0\}$.
The coefficient $V_p^{p-1}$ belongs to $C^{0,\alpha}$ because
$V_p$ is Lipschitz and $\alpha<p-1$. The regular
Dirichlet ODE and its Schauder estimate therefore give, for each
$f\in C^{0,\alpha}$, a unique zero-trace solution of
$-A_Dw=f$ and the bound
\[
 \|w\|_{C^{2,\alpha}}\le C\|f\|_{C^{0,\alpha}}.
\]
In particular, $D_v\mathcal E(0,0,V_p)$ is an isomorphism.
The implicit function theorem yields a $C^1$ map
$(\delta,b)\mapsto v(\delta,b)$ near $(0,0)$, with
$v(0,0)=V_p$, and gives local uniqueness of solutions near this
profile. We now restrict this map to the physical trace by setting
\[
 V_\delta=v(\delta,\delta^a)+\delta^a(1-t).
\]
The first derivatives of $v$ are bounded on a smaller parameter
neighborhood. Integrating along the segment from $(0,0)$ to
$(\delta,\delta^a)$, we obtain
\[
 \|v(\delta,\delta^a)-V_p\|_{C^{2,\alpha}}
 \le C(\delta+\delta^a).
\]
Adding the boundary lift proves \eqref{eq:V-estimate}. The
local uniqueness supplied by the theorem applies to any solution of
\eqref{eq:Vdelta} sufficiently close to $V_p$, after subtracting
the same lift $\delta^a(1-t)$.

We next check that the extension by the positive part has introduced
no change in the desired equation. By Lemma~\ref{lem:pulse},
$V_p'(0)>0$ and $V_p'(1)<0$. Choose
$\epsilon\in(0,1/2)$ and $c>0$ so that
$V_p'\ge2c$ on $[0,\epsilon]$ and
$V_p'\le-2c$ on $[1-\epsilon,1]$.
The $C^1$ convergence in \eqref{eq:V-estimate} gives the
corresponding bounds $V_\delta'\ge c$ and
$V_\delta'\le-c$ for small $\delta$. Writing $b=\delta^a$,
we therefore have
\[
 V_\delta(t)\ge b+ct\quad(0\le t\le\epsilon),\qquad
 V_\delta(t)\ge c(1-t)\quad(1-\epsilon\le t\le1).
\]
On $[\epsilon,1-\epsilon]$, uniform convergence to the positive
function $V_p$ gives a positive lower bound independent of small
$\delta$. Thus $V_\delta>0$ on $[0,1)$, its zero at $1$
is simple, and $(V_\delta)_+^p=V_\delta^p$. This proves that
the profile solves \eqref{eq:Vdelta} with the original power.

Although the curve $b=\delta^a$ need not be differentiable at
zero, it is smooth for $\delta>0$. The two-parameter representation
gives, for these positive parameters,
\[
 \delta\partial_\delta V_\delta
 =\delta\,\partial_\delta v(\delta,\delta^a)
   +a\delta^a\bigl(\partial_bv(\delta,\delta^a)+(1-t)\bigr).
\]
The bounded parameter derivatives prove
\eqref{eq:V-derivative}.

To obtain smooth dependence for each fixed small $\delta_0>0$,
we use the simple outer zero more precisely. The quotient
$V_{\delta_0}/(1-t)$ extends to a strictly positive
$C^{1,\alpha}$ function on $[0,1]$. For all nearby profiles
with zero right trace, division by $1-t$ is bounded into
$C^{1,\alpha}$, as in \eqref{eq:division}. Their quotients
remain bounded away from zero. Writing the power as
$(1-t)^p(V/(1-t))^p$ therefore makes the residual $C^\infty$
on this neighborhood, also when $p$ is not an integer.
Moreover, the coefficients $V_\delta^{p-1}$ converge to
$V_p^{p-1}$ in $C^{0,\alpha}$: they converge uniformly and
are uniformly bounded in $C^{0,\beta}$, so the same interpolation
argument applies. Together with
$\delta/(\lambda+\delta t)\to0$ and $\delta^2\to0$, this
shows that $D_v\mathcal E(\delta,\delta^a,v(\delta,\delta^a))$
converges in operator norm to the above isomorphism. It is thus
invertible for all sufficiently small $\delta$. The smooth implicit
function theorem near $(\delta_0,\delta_0^a)$ gives a smooth
solution map there. By local uniqueness, this is the already
constructed family, which is consequently $C^\infty$ on the open
interval of small positive widths.
\end{proof}

Returning to the original variables, we define
\begin{equation}\label{eq:Udelta}
 U_\delta(r)=\delta^{-a}
       V_\delta\left(\frac{r-\lambda}{\delta}\right).
\end{equation}
The scaling calculation and positivity above show that this function
solves \eqref{eq:radial-original}, is positive in the annulus, and
has inner value $1$. Its outer normal derivative is
\begin{equation}\label{eq:ddelta}
 d_\delta
 =U_\delta'(\lambda+\delta)
 =\delta^{-a-1}V_\delta'(1)<0.
\end{equation}
Thus the radial Bernoulli constant is $d_\delta^2>0$. At the
inner circle, $C^1$ convergence gives $U_\delta'(\lambda)>0$,
and uniform convergence of the profiles gives
\[
 \max U_\delta=\delta^{-a}(M_p+o(1)).
\]
In particular, these radial solutions lie outside the ordered class
$0<U_\delta<1$.

\subsubsection{A derivative estimate for the linearized potential}
To verify transversality, we also need an estimate for the parameter
derivative of the potential. When $1<p<2$, this derivative contains
the factor $V_\delta^{p-2}$, which is singular at the outer zero.
The profile variation vanishes at that endpoint, and its value at
the inner endpoint is controlled by the small rescaled trace. We
use both facts in the following estimate.
Now let us establish convergence of the potential and its parameter derivative.

\begin{lemma}
\label{lem:potential-derivative}
As $\delta\downarrow0$,
\begin{equation}\label{eq:potential-convergence}
 \|V_\delta^{p-1}-V_p^{p-1}\|_{L^\infty}\longrightarrow0,\qquad
 \|\delta\partial_\delta(V_\delta^{p-1})\|_{L^\infty}
 \longrightarrow0.
\end{equation}
Moreover, the second norm is at most $C(\delta+\delta^a)$.
\end{lemma}

\begin{proof}
The profiles are uniformly bounded and converge uniformly to $V_p$.
Uniform continuity of $s\mapsto s^{p-1}$ on a bounded nonnegative
interval therefore proves the first limit. Set
$Z_\delta=\delta\partial_\delta V_\delta$ and $b=\delta^a$.
Lemma~\ref{lem:radial-family} and the differentiated traces give
\[
 \|Z_\delta\|_{C^1}\le C(\delta+b),\qquad
 Z_\delta(0)=ab,\quad Z_\delta(1)=0.
\]
For $0<t<1$, the ordinary chain rule gives
\[
 \delta\partial_\delta(V_\delta^{p-1})
 =(p-1)V_\delta^{p-2}Z_\delta.
\]
If $p\ge2$, the profiles and the factor $V_\delta^{p-2}$
are bounded, so \eqref{eq:V-derivative} proves the required
bound directly.

Suppose now that $1<p<2$. The positivity proof gives fixed
$c,\epsilon>0$ for which
\[
 V_\delta(t)\ge c(t+b)\quad(0\le t\le\epsilon),\qquad
 V_\delta(t)\ge c(1-t)\quad(1-\epsilon\le t\le1),
\]
and gives a uniform positive lower bound on the middle interval.
Since $Z_\delta(1)=0$, integrating its first derivative from
the right endpoint yields
$|Z_\delta(t)|\le C(\delta+b)(1-t)$. Consequently,
\[
 V_\delta(t)^{p-2}|Z_\delta(t)|
 \le C(\delta+b)(1-t)^{p-1}
 \le C(\delta+b)
 \quad(1-\epsilon\le t<1).
\]
In particular, this product extends continuously by zero at $t=1$.
At the left endpoint, we instead integrate from $0$ and use
$Z_\delta(0)=ab$, obtaining
$|Z_\delta(t)|\le Cb+C(\delta+b)t$. Thus
\[
 \begin{split}
 V_\delta(t)^{p-2}|Z_\delta(t)|
 &\le C\frac{b+(\delta+b)t}{(t+b)^{2-p}}\\
 &\le Cb^{p-1}+C(\delta+b)t^{p-1}
 \le C\bigl(b^{p-1}+\delta+b\bigr)
 \quad(0\le t\le\epsilon).
 \end{split}
\]
The second inequality follows from $t+b\ge b$ in the first
term and $t+b\ge t$ in the second. Since
$b^{p-1}=\delta^{a(p-1)}=\delta^2$, this is bounded by
$C(\delta+b)$ for $\delta\le1$. The uniform positive lower
bound for $V_\delta$ and the $C^0$ estimate on $Z_\delta$
give the same bound on the middle interval. After multiplying by
$p-1$, we obtain the stated norm estimate on the whole interval.
Its right-hand side tends to zero, which proves the second limit.
\end{proof}

\subsubsection{Radial Dirichlet invertibility}
The radial Dirichlet linearization in the rescaled variables is
\begin{equation}\label{eq:radial-Jacobi-scaled}
 A_{\delta,0}h
 =-h''-\frac{\delta}{\lambda+\delta t}h'
   +\bigl(\delta^2-pV_\delta^{p-1}\bigr)h,\qquad h(0)=h(1)=0.
\end{equation}
Put $w_\delta(t)=1+\delta t/\lambda$. Since
$w_\delta'/w_\delta=\delta/(\lambda+\delta t)$, we can
write the derivative part as
$-w_\delta^{-1}(w_\delta h')'$. Integration by parts with
the zero endpoint values then gives the quadratic form
\[
 q_\delta(h)
 =\int_0^1w_\delta\bigl(|h'|^2+
       (\delta^2-pV_\delta^{p-1})h^2\bigr)\dd t,
 \qquad h\in H_0^1(0,1),
\]
with squared norm $\int_0^1w_\delta h^2\dd t$.
Thus $A_{\delta,0}$ is self-adjoint in this weighted
$L^2$ space and has a discrete sequence of real eigenvalues.

The weight tends uniformly to $1$, and the potential tends uniformly
to $-pV_p^{p-1}$ by Lemma~\ref{lem:potential-derivative}.
If $q_D(h)=\int_0^1(|h'|^2-pV_p^{p-1}h^2)\dd t$,
these convergences give
\[
 |q_\delta(h)-q_D(h)|
 \le\epsilon_\delta\bigl(\|h'\|_{L^2}^2+
                              \|h\|_{L^2}^2\bigr),\qquad
 \epsilon_\delta\longrightarrow0,
\]
as well as
$\int_0^1w_\delta h^2\dd t=(1+O(\delta))\|h\|_{L^2}^2$,
uniformly for $h$. All potentials are uniformly bounded. We may
therefore add a fixed sufficiently large multiple of the squared
norm to the forms and compare the resulting positive Rayleigh
quotients. The min--max principle shows that each fixed eigenvalue
of $A_{\delta,0}$ converges to the corresponding eigenvalue of
$A_D$. By Lemma~\ref{lem:pulse-D}, the first limiting
eigenvalue is negative and the second is strictly positive.
For sufficiently small $\delta$, the same signs hold for the first
two eigenvalues of $A_{\delta,0}$. Every higher eigenvalue is at
least the second one. Hence $A_{\delta,0}$ has exactly one
negative eigenvalue and no zero eigenvalue.

We can also see the absence of a zero eigenvalue directly. If
$A_{\delta_n,0}h_n=0$ with $\delta_n\downarrow0$, normalize
the solutions by $\|h_n\|_{C^0}=1$. The uniformly bounded
potentials, together with the Dirichlet equation, give a uniform
$C^2$ bound; for example, the one-dimensional elliptic estimate
$\|h_n\|_{W^{2,q}}\le C(\|A_{\delta_n,0}h_n\|_{L^q}
+\|h_n\|_{L^q})$ first bounds $h_n'$ for any $q>1$, and
the equation then bounds $h_n''$. A subsequence converges in
$C^1$ to a zero-trace function $h$ with $\|h\|_{C^0}=1$.
Passing to the equation using the uniform coefficient convergence
gives $A_Dh=0$, contrary to Lemma~\ref{lem:pulse-D}.
This also explains why radial Dirichlet invertibility persists along
the large annular family.

The preceding construction reduces the thin-annulus problem to a
fixed interval. We now prove directly the one-dimensional profile
properties used above and the mixed--Dirichlet gap that will determine
the critical angular scale.

\subsection{The limiting profile for the radial branch}\label{sec:pulse}

The rescaling in the preceding subsection explains the form of the
limiting profile. In \eqref{eq:Vdelta}, the coefficients of the
curvature term and the linear reaction term tend to zero as
$\delta\to0$. At the same time, the rescaled inner value is
$\delta^a\to0$, although the original inner value is $1$. The
choice $a(p-1)=2$ keeps the second derivative and the power term
at the same scale. Thus the nonzero limit of the large radial branch
solves
\begin{equation}\label{eq:pulse-intro}
 -V''=V^p\quad\text{in }(0,1),\qquad V(0)=V(1)=0.
\end{equation}
We now prove the properties of this profile used in the radial
construction. We also compare its Dirichlet and mixed spectra, since
the state equation and the linearized outer boundary condition impose
different conditions at $t=1$.

\subsubsection{Existence, uniqueness, and elementary identities}
The following lemma characterizes the positive limiting profile.

\begin{lemma}\label{lem:pulse}
For each $p>1$, problem \eqref{eq:pulse-intro} admits a unique
positive solution $V=V_p$. This solution is symmetric about
$t=1/2$, strictly concave, strictly increasing on $(0,1/2)$,
and strictly decreasing on $(1/2,1)$. Moreover,
$V'(0)>0$ and $V'(1)=-V'(0)<0$.
\end{lemma}

\begin{proof}
We first construct a positive trajectory with prescribed maximum
$M>0$. We solve the initial-value problem with value $M$ and
derivative zero at its maximum. For ODE existence and uniqueness,
we may extend the nonlinearity to $(v_+)^p$ on the real line;
this extension is $C^1$ since $p>1$. While the trajectory is
positive, it solves the original equation and has strictly negative
second derivative. Multiplication of the equation by $V'$ gives
\begin{equation}\label{eq:pulse-energy}
 \frac12V'(t)^2+\frac1{p+1}V(t)^{p+1}
 =\frac1{p+1}M^{p+1}.
\end{equation}
On its decreasing part, the derivative is therefore
\[
 V'=-\sqrt{\frac{2}{p+1}}
       \sqrt{M^{p+1}-V^{p+1}}.
\]
Separation of variables shows that the time from the maximum to
the first zero is
\[
 T(M)=\sqrt{\frac{p+1}{2}}
       \int_0^M\frac{\dd v}{\sqrt{M^{p+1}-v^{p+1}}}
 =\sqrt{\frac{p+1}{2}}\,M^{-(p-1)/2}
       \int_0^1\frac{\dd z}{\sqrt{1-z^{p+1}}}.
\]
This integral is finite: its only singularity is at $z=1$, where
the integrand is bounded by a constant times $(1-z)^{-1/2}$.
In particular, the trajectory reaches zero in finite time and its
derivative there is nonzero. Reflection about the maximum gives
the increasing part, with the same length $T(M)$, by uniqueness
of the autonomous initial-value problem.

The total distance between the two zeros is $2T(M)$. It is a
continuous, strictly decreasing function of $M$, tends to infinity
as $M\downarrow0$, and tends to zero as $M\to\infty$. There is
therefore exactly one $M$ for which $2T(M)=1$. Translating this
trajectory to put its first zero at $0$ gives a positive solution
of \eqref{eq:pulse-intro}, with its maximum at $1/2$.

Conversely, any positive solution of \eqref{eq:pulse-intro} has
$V''=-V^p<0$ in the interior. Its derivative is strictly
decreasing and its maximum is attained at a unique interior point.
The two endpoint derivatives cannot vanish: zero value and zero
derivative at either endpoint would give the identically zero
solution by ODE uniqueness. They consequently have the signs
$V'(0)>0$ and $V'(1)<0$. Applying \eqref{eq:pulse-energy} to
this solution shows that the distances from its maximum to both
endpoints are $T(M)$. Hence $2T(M)=1$, and its maximum and
initial slope agree with those of the constructed solution. ODE
uniqueness proves uniqueness of the positive solution. Reflection
about $1/2$ proves symmetry, and \eqref{eq:pulse-energy} gives
$V'(1)=-V'(0)$. Strict concavity then gives the asserted
monotonicity on the two half-intervals.
\end{proof}

Let $I_p=\int_0^1(1-z^{p+1})^{-1/2}\dd z$. Solving
$2T(M_p)=1$ and evaluating \eqref{eq:pulse-energy} at $t=0$
give
\begin{equation}\label{eq:pulse-constants}
 M_p=\bigl(\sqrt{2(p+1)}\,I_p\bigr)^{2/(p-1)},\qquad
 V_p'(0)=\sqrt{\frac{2}{p+1}}\,M_p^{(p+1)/2}.
\end{equation}
These formulas determine the leading size of the radial solution
and its Bernoulli constant. The simple endpoint zeros also show that
$V_p(t)$ is comparable to $t$ near $0$ and to $1-t$ near $1$.
The constant $\kappa_p$ governing the critical radii is determined
by the spectral problem \eqref{eq:DN-intro}.

\subsubsection{Dirichlet nondegeneracy}
Consider the self-adjoint Dirichlet operator
\[
 A_D=-\frac{d^2}{dt^2}-pV_p^{p-1}
 \quad\text{on }H^2(0,1)\cap H_0^1(0,1).
\]
Its potential is continuous on $[0,1]$, so this is a regular
Sturm--Liouville operator. Let $\nu_D$ denote its lowest eigenvalue.
We next describe the Dirichlet spectrum of the limiting profile.

\begin{lemma}\label{lem:pulse-D}
The operator $A_D$ has exactly one negative eigenvalue and no zero
eigenvalue. In particular, $\nu_D<0$.
\end{lemma}

\begin{proof}
We first test the Dirichlet quadratic form with $V=V_p$.
Multiplying $-V''=V^p$ by $V$ and integrating by parts, with
zero boundary values, gives $\int_0^1V'^2\dd t
=\int_0^1V^{p+1}\dd t$. Therefore
\[
 \int_0^1\bigl(V'^2-pV^{p+1}\bigr)\dd t
 =(1-p)\int_0^1V^{p+1}\dd t<0.
\]
The variational characterization of $\nu_D$ gives $\nu_D<0$.

To exclude zero and count the negative eigenvalues, we use the
scaling of the profile equation. Locally in the scaling parameter
$s>0$, the function $s^aV(st)$ solves the same differential
equation, wherever its argument lies in $(0,1)$, because
$a+2=ap$. Its derivative at $s=1$ is
\begin{equation}\label{eq:scaling-Jacobi}
 Z(t)=aV(t)+tV'(t),\qquad a=\frac{2}{p-1}.
\end{equation}
We can also verify its equation directly. In the interior,
\[
 \begin{aligned}
 Z'&=(a+1)V'+tV'',\\
 Z''&=(a+2)V''+tV'''
       =-(a+2)V^p-ptV^{p-1}V',\\
 Z''+pV^{p-1}Z&=(pa-a-2)V^p=0.
 \end{aligned}
\]
The equation extends to the endpoints by continuity. Moreover,
\[
 Z(0)=0,\qquad Z'(0)=(a+1)V'(0)>0,\qquad Z(1)=V'(1)<0.
\]
Every solution of the zero-energy equation with zero left endpoint
value is determined by its left derivative. It is consequently a
multiple of $Z$. Since $Z(1)\ne0$, no nonzero such solution
satisfies the right Dirichlet condition, and zero is not an
eigenvalue of $A_D$.

It remains to count the zeros of $Z$. On $(0,1/2]$, we have
$V>0$ and $V'\ge0$, so $Z>0$. On $(1/2,1)$,
\[
 Z'=(a+1)V'+tV''<0,
\]
because both $V'$ and $V''$ are negative. Since
$Z(1/2)=aV(1/2)>0$ and $Z(1)<0$, the function $Z$ has
exactly one interior zero. The Sturm oscillation theorem
\cite{Zettl}, applied to the solution of $A_DZ=0$ with zero
left value and positive left derivative, identifies the number of
interior zeros with the number of Dirichlet eigenvalues below zero;
here the right endpoint is not a zero. Thus $A_D$ has exactly
one negative eigenvalue.
\end{proof}

In particular, the homogeneous Dirichlet linearization at $V_p$
has trivial kernel. Since the potential belongs to $C^{0,\alpha}$
under \eqref{eq:alpha}, the Dirichlet Fredholm alternative and
the one-dimensional Schauder estimate give the isomorphism used
in the preceding implicit-function construction of the radial family.

\subsubsection{The mixed eigenvalue and the strict gap}
On the form domain
$E=\{h\in H^1(0,1):h(0)=0\}$, define
\[
 q_{\mathrm{DN}}(h)=
 \int_0^1\bigl(|h'|^2-pV_p^{p-1}h^2\bigr)\dd t.
\]
The associated lowest eigenvalue is the infimum of
$q_{\mathrm{DN}}(h)/\|h\|_{L^2}^2$ over nonzero $h\in E$.
We denote its operator by $A_{\mathrm{DN}}$; it has the same
differential expression as $A_D$, but its domain is
\[
 \{h\in H^2(0,1):h(0)=0,\ h'(1)=0\}.
\]
Thus its eigenvalue problem is \eqref{eq:DN-intro}.
The following lemma establishes a strict separation between the principal mixed and Dirichlet eigenvalues.

\begin{lemma}\label{lem:gap}
The principal mixed eigenvalue is simple, has a positive
eigenfunction, and satisfies
\begin{equation}\label{eq:gap}
 \nu_{\mathrm{DN}}<\nu_D<0,\qquad
 g_p:=\nu_D-\nu_{\mathrm{DN}}>0.
\end{equation}
If $\phi_p$ is its positive eigenfunction normalized by
$\phi_p(1)=1$, then
\[
 \kappa_p^2=-\nu_{\mathrm{DN}},\qquad
 A_D+\kappa_p^2\ge g_p\Id
\]
in the Dirichlet quadratic-form sense.
\end{lemma}

\begin{proof}
Put $q(t)=pV_p(t)^{p-1}$, a bounded continuous function. We
minimize the mixed Rayleigh quotient under the normalization
$\|h\|_{L^2}=1$. Since
\[
 q_{\mathrm{DN}}(h)\ge\|h'\|_{L^2}^2-\|q\|_{L^\infty}
 \qquad\text{when }\|h\|_{L^2}=1,
\]
a minimizing sequence is bounded in $H^1(0,1)$. We pass to a
subsequence converging weakly in $H^1$ and strongly in $L^2$.
The limit belongs to $E$ and has $L^2$ norm one. Weak lower
semicontinuity of the derivative term and strong convergence of
the potential term show that it attains the infimum.

We may replace this minimizer by its absolute value, which preserves
both its norm and its energy. The resulting nonnegative minimizer
$\phi$ satisfies the Euler--Lagrange identity
\[
 \int_0^1(\phi'h'-q\phi h)\dd t
 =\nu_{\mathrm{DN}}\int_0^1\phi h\dd t
 \qquad(h\in E).
\]
Interior variations give $-\phi''-q\phi=\nu_{\mathrm{DN}}\phi$.
ODE regularity gives $\phi\in C^2([0,1])$, and integration
by parts in the displayed identity, with arbitrary $h(1)$,
gives the natural condition $\phi'(1)=0$. A zero of $\phi$
in the interior would be a minimum and would also satisfy
$\phi'=0$ there. ODE uniqueness would then give $\phi\equiv0$,
contrary to its normalization. Likewise, $\phi(1)=0$ together
with $\phi'(1)=0$ is impossible. Thus $\phi>0$ on $(0,1]$,
and we may normalize it by $\phi(1)=1$. Any eigenfunction for
this eigenvalue has zero value at $0$ and is determined by its
derivative there. ODE uniqueness therefore makes it a multiple
of $\phi$, proving simplicity.

The Dirichlet form domain $H_0^1(0,1)$ is contained in $E$,
and the two forms have the same expression. Hence
$\nu_{\mathrm{DN}}\le\nu_D$. To prove strictness, suppose
that equality holds and choose a normalized principal Dirichlet
eigenfunction $\psi$. Its Rayleigh quotient is then the mixed
infimum, so it also satisfies the mixed Euler--Lagrange identity
for every $h\in E$. Integration by parts yields $\psi'(1)=0$.
Its Dirichlet condition already gives $\psi(1)=0$, so ODE
uniqueness again forces $\psi\equiv0$, a contradiction.
Together with Lemma~\ref{lem:pulse-D}, this proves
\eqref{eq:gap} and makes $\kappa_p$ well defined.

Finally, for every $h\in H_0^1(0,1)$, the Dirichlet variational
inequality gives
\[
 \int_0^1\bigl(|h'|^2-pV_p^{p-1}h^2+\kappa_p^2h^2\bigr)\dd t
 \ge(\nu_D-\nu_{\mathrm{DN}})\|h\|_{L^2}^2
 =g_p\|h\|_{L^2}^2.
\]
This is the asserted strictly positive Dirichlet form bound.
\end{proof}

The gap $g_p$ explains how the outer boundary equation can have
a kernel while the Dirichlet state equation remains invertible.
At the critical angular scale, the limiting mixed operator is
$A_{\mathrm{DN}}+\kappa_p^2$, with principal eigenvalue zero,
whereas $A_D+\kappa_p^2$ is bounded below by $g_p\Id$.
Together with the radial nondegeneracy above, this separation
provides the limiting spectral bounds used later to exclude
Dirichlet resonance at the bifurcation radii.

\subsection{The state equation and the Bernoulli equation}
For the remainder of the construction, the trivial state is
$U_{R-\lambda}$ from \eqref{eq:Udelta}. We now set up the equations
near one such state on $A_*$; the spectral analysis below will select
$R_*=R_K$ and verify the required Dirichlet invertibility.

Choose a smooth radial cutoff $\chi$, equal to zero near $\lambda$
and to one near $R_*$. For $(R,\eta)$ near $(R_*,0)$, set
\begin{equation}\label{eq:power-pullback}
 q_{R,\eta}(r,\theta)=r+\chi(r)(R-R_*+\eta(\theta)),\qquad
 T_{R,\eta}(r,\theta)=(q_{R,\eta}(r,\theta),\theta).
\end{equation}
The smallness condition makes $q_r>0$ and
$R+\eta>\lambda$, so $T_{R,\eta}$ maps $A_*$ diffeomorphically
onto $\Omega_{R,\eta}\setminus\overline{B_\lambda}$. It is the
identity in an inner collar. To make the pulled-back equation
explicit, write $q=q_{R,\eta}$ and
\[
 \mathcal D_q=\partial_\theta-\frac{q_\theta}{q_r}\partial_r.
\]
The chain rule in polar coordinates gives
\begin{equation}\label{eq:power-pulled-laplacian}
 \Delta^p_{R,\eta}v
 =\frac1{q_r}\partial_r\left(\frac{v_r}{q_r}\right)
    +\frac{v_r}{q q_r}+\frac1{q^2}\mathcal D_q^2v.
\end{equation}
Thus the Dirichlet equation for this specific power law is
\begin{equation}\label{eq:power-state-equation}
 \mathcal E_p(R,\eta,w)
 =-\Delta^p_{R,\eta}(b_*+w)+(b_*+w)-(b_*+w)^p=0,
\end{equation}
where $w\in Z_{p,K}$ and $b_*+w$ lies in a sufficiently small
$C^{2,\alpha}$ neighborhood of $U_{R_*-\lambda}$.
This residual takes values in $W_{p,K}$. Formula
\eqref{eq:power-pulled-laplacian}, multiplication in H\"older spaces,
and Lemma~\ref{lem:power-map} show that $\mathcal E_p$ is
$C^\infty$. Its coefficients depend smoothly on $R,\eta$ and
their indicated derivatives, and uniform ellipticity follows from
$q_r>0$ and proximity to the identity.

At the trivial state, $w_*=U_{R_*-\lambda}-b_*$, the derivative
with respect to the state is the concrete Dirichlet operator
\begin{equation}\label{eq:power-state-linearization}
 D_w\mathcal E_p(R_*,0,w_*)w
 =-\Delta w+\bigl(1-pU_{R_*-\lambda}^{p-1}\bigr)w,
 \qquad w\in Z_{p,K}.
\end{equation}
The potential belongs to $C^{0,\alpha}$ by \eqref{eq:alpha}.
The Dirichlet Laplacian is an isomorphism
$Z_{p,K}\to W_{p,K}$ by the Schauder theory \cite{GT} and
uniqueness, and multiplication by this potential is a compact
perturbation, using the compact embedding
$C^{2,\alpha}\hookrightarrow C^{0,\alpha}$. Consequently, the operator defined in
\eqref{eq:power-state-linearization} is Fredholm of index zero.
We next construct the state map for the power law.

\begin{proposition}\label{prop:power-state}
Suppose that the Dirichlet operator in
\eqref{eq:power-state-linearization} has trivial kernel in $Z_{p,K}$.
Then \eqref{eq:power-state-equation} defines a unique nearby state
\[
 (R,\eta)\longmapsto v_{R,\eta}=b_*+w(R,\eta)
\]
which is positive in $A_*$ and $C^\infty$ in
$C^{2,\alpha}(\overline A_*)$.
For $\eta=0$, it equals $U_{R-\lambda}\circ T_{R,0}$.
\end{proposition}

\begin{proof}
The Fredholm statement and the kernel assumption make
\eqref{eq:power-state-linearization} a bounded isomorphism. The
implicit function theorem applied to \eqref{eq:power-state-equation}
gives the state map and uniqueness. The neighborhood is chosen small enough that Lemma~\ref{lem:power-map}
guarantees positivity and a negative outer normal derivative. The constructed radial
family solves the same equation when $\eta=0$; local uniqueness
identifies the two maps there.
\end{proof}

Set $u_{R,\eta}=v_{R,\eta}\circ T_{R,\eta}^{-1}$. Its outer
normal derivative is particularly simple for
\eqref{eq:power-pullback}. Since $\chi=1$ near $R_*$ and the
outer Dirichlet trace is zero, the chain rule and the outer normal give
\begin{equation}\label{eq:power-normal}
 N_p(R,\eta)(\theta)
 =\sqrt{1+\frac{\eta_\theta(\theta)^2}{(R+\eta(\theta))^2}}
                 \,\partial_r v_{R,\eta}(R_*,\theta).
\end{equation}
Indeed, $q_r=1$ on the reference outer circle, while differentiation
of the zero physical trace gives
$u_\theta(R+\eta,\theta)=-\eta_\theta u_r(R+\eta,\theta)$.
We define the boundary equation by
\begin{equation}\label{eq:power-Bernoulli}
 F_{p,K}(R,\eta)=\Pi\bigl(N_p(R,\eta)^2\bigr)=0,
 \qquad \Pi g=g-\av g.
\end{equation}
The trace estimate and \eqref{eq:power-normal} imply that
$F_{p,K}$ is a $C^\infty$ map from a neighborhood of $(R_*,0)$
in $\R\times X_K$ to $Y_K$, and
$F_{p,K}(R,0)=0$ on the radial branch. Since the boundary trace
is constant, its tangential derivative vanishes, and the boundary
condition in \eqref{eq:power} is exactly
\eqref{eq:power-Bernoulli}, with
\begin{equation}\label{eq:power-recover-Q}
 Q=\av{N_p(R,\eta)^2}>0.
\end{equation}
The strict inequality follows from the negative outer derivative at
the trivial state and continuity. The map is equivariant under
reflection and rotation, including the half-period rotation
$\eta(\theta)\mapsto\eta(\theta+\pi/K)$ on $X_K$.

\subsection{Linearization, the boundary symbol, and the range}
We set $A_{\lambda,R}=B_R\setminus\overline{B_\lambda}$.
Write $U=U_{R-\lambda}$ and $d=U'(R)<0$. In this subsection,
assume that the power-law Dirichlet linearization is invertible in
the $D_K$ class. We will verify this at the radii used for bifurcation.
The radial equation and $U(R)=0$ give
\begin{equation}\label{eq:outer-second}
 U''(R)=-\frac dR.
\end{equation}
For $\eta=\epsilon\xi$, let $\dot v$ be the material derivative.
At a reference annulus of radius $R$, the Eulerian derivative is
$u'=\dot v-\chi(r)U'(r)\xi(\theta)$. Differentiating the state
equation and its boundary traces gives
\begin{equation}\label{eq:shape-state}
 \begin{cases}
 (-\Delta+1-pU^{p-1})u'=0&\text{in }A_{\lambda,R},\\
 u'=0&\text{on }r=\lambda,\\
 u'=-d\xi&\text{on }r=R.
 \end{cases}
\end{equation}
The material derivative has zero boundary traces; the nonzero outer
trace of $u'$ accounts for the boundary displacement. Differentiating
\eqref{eq:power-normal}, or differentiating the physical normal,
yields
\begin{equation}\label{eq:shape-normal}
 D_\eta N_p(R,0)[\xi]
       =\partial_r u'(R,\theta)+U''(R)\xi(\theta).
\end{equation}
The normal-vector variation is tangential and has zero scalar product
with the radial gradient at the trivial state.

For $j=nK$, $n\ge1$, let $H_j=H_j(R;\cdot)$ solve
\begin{equation}\label{eq:Hj}
 -H_j''-\frac1rH_j'
 +\left(\frac{j^2}{r^2}+1-pU^{p-1}\right)H_j=0,
 \qquad H_j(\lambda)=0,\quad H_j(R)=1.
\end{equation}
Dirichlet invertibility gives existence and uniqueness. Substitution
of $\xi=\cos(j\theta)$ into \eqref{eq:shape-state} gives
$u'=-dH_j(r)\cos(j\theta)$. Define
\begin{equation}\label{eq:Tj}
 T_j(R)=H_j'(R)+\frac1R.
\end{equation}
Equations \eqref{eq:outer-second}--\eqref{eq:shape-normal} then give
\begin{equation}\label{eq:multiplier}
 D_\eta F_{p,K}(R,0)[\cos(j\theta)]
      =m_j(R)\cos(j\theta),\qquad m_j(R)=-2d^2T_j(R).
\end{equation}
In particular, $T_j(R)=0$ precisely when the concrete mixed problem
\begin{equation}\label{eq:mixed-general}
 \left[-\frac{d^2}{dr^2}-\frac1r\frac d{dr}
       +\frac{j^2}{r^2}+1-pU^{p-1}\right]h=0,
 \qquad h(\lambda)=0,\quad h'(R)+R^{-1}h(R)=0
\end{equation}
has a nonzero solution. The Robin coefficient is produced by the
radial equation at the outer boundary.

Let $\Lambda_p$ be the outer Dirichlet-to-Neumann map for
$-\Delta+1-pU^{p-1}$ with zero inner trace. At this radial state,
\[
 L_{p,K}:=D_\eta F_{p,K}(R,0)
       =-2d^2(\Lambda_p+R^{-1}\Id):X_K\longrightarrow Y_K.
\]
The following lemma establishes the Fredholm property and characterizes the range for the power law.

\begin{lemma}
\label{lem:power-Fredholm}
The operator $L_{p,K}:X_K\to Y_K$ is Fredholm of index zero.
If $T_K(R)=0$ and $T_{nK}(R)\ne0$ for all $n\ge2$, then
\begin{align*}
 \Ker L_{p,K}&=\Span\{\cos(K\theta)\},\\
 \Ran L_{p,K}&=\left\{g\in Y_K:
       \int_0^{2\pi}g(\theta)\cos(K\theta)\dd\theta=0\right\}.
\end{align*}
\end{lemma}

\begin{proof}
We write $c(r)=1-pU(r)^{p-1}$ and $\mathcal L_p=-\Delta+c$.
Our choice of $\alpha<p-1$ gives
$c\in C^{0,\alpha}(\overline A_{\lambda,R})$. The assumed
Dirichlet invertibility and the Schauder estimates define the
unique extensions $P_p\xi$ by
\[
 \mathcal L_p(P_p\xi)=0\quad\text{in }A_{\lambda,R},\qquad
 P_p\xi|_{r=\lambda}=0,\quad P_p\xi|_{r=R}=\xi.
\]
Uniqueness preserves the $D_K$ symmetries and makes their angular
means zero. Thus $\Lambda_p\xi=\partial_r(P_p\xi)(R,\cdot)$
is bounded from $X_K$ to $Y_K$.

Let $J:X_K\hookrightarrow Y_K$ be the natural inclusion.
We apply the Dirichlet--Robin argument in
Lemma~\ref{lem:Fredholm} with potential $c$ and outer condition
$\partial_r v+R^{-1}v=g$. It uses only the stated H\"older
regularity of $c$. Since $dS=R\dd\theta$ on the outer circle,
the Fredholm orthogonality conditions are unchanged by using
$\dd\theta$. Therefore $\Lambda_p+R^{-1}J$ has closed range
equal to the $L^2(\dd\theta)$ orthogonal complement of its
finite-dimensional kernel within $Y_K$, and its index is zero. The nonzero factor
$-2d^2$ gives the same conclusion for $L_{p,K}$.

Taking angular Fourier coefficients in the extension equation
and using Dirichlet uniqueness gives, for $j=nK$,
\[
 \frac1\pi\int_0^{2\pi}(L_{p,K}\eta)(\theta)
                            \cos(j\theta)\dd\theta
 =m_j(R)\frac1\pi\int_0^{2\pi}\eta(\theta)
                            \cos(j\theta)\dd\theta.
\]
Under the stated hypotheses, only the $K$-th multiplier vanishes.
The symmetries, zero mean, and uniqueness of Fourier coefficients
therefore give $\Ker L_{p,K}=\Span\{\cos(K\theta)\}$.
The same identity shows that the range lies in the stated
hyperplane. Both have codimension one by the Fredholm index, so
they coincide. This argument does not require convergence of
Fourier partial sums in the H\"older norm.
\end{proof}

For later use, the principal eigenvalue $\zeta_j(R)$ of
\eqref{eq:mixed-general} has the Rayleigh quotient
\begin{equation}\label{eq:mixed-form}
 \zeta_j(R)=\min_{\substack{h\in H^1(\lambda,R)\setminus\{0\}\\
                              h(\lambda)=0}}
 \frac{\displaystyle\int_\lambda^R
       \left(r|h'|^2+\frac{j^2}{r}h^2+r(1-pU^{p-1})h^2\right)\dd r
       +h(R)^2}
      {\displaystyle\int_\lambda^Rrh^2\dd r}.
\end{equation}
The positive endpoint term gives the Robin coefficient $R^{-1}$.
Comparing the quotients at the same state and radius gives
\begin{equation}\label{eq:mode-order}
 \zeta_j(R)\ge\zeta_\ell(R)+\frac{j^2-\ell^2}{R^2}
                 >\zeta_\ell(R),\qquad j>\ell.
\end{equation}
Thus a zero of the principal $K$-eigenvalue isolates the boundary
kernel in the $D_K$ class. It remains to construct such a zero,
prove Dirichlet nonresonance there, and check that it is transverse
as $R$ varies along the radial branch.

\subsection{Verification of the spectral and transversality conditions}\label{sec:crossing}

We now construct the critical radii for the power law. We fix
the angular frequency $K$ and vary the outer radius in the regime
where $K\delta/\lambda$ remains of order one. We first locate a
simple zero of the principal mixed eigenvalue. We then verify that
the Dirichlet state problem is invertible at this radius and that
the zero gives the required transverse crossing of the Bernoulli
operator. All derivatives with respect to the radius below are
taken along the radial family $U_{R-\lambda}$.

\subsubsection{The scaled mixed operator}
For the radial solution $U_\delta$, the Jacobi operator in the
$j$-th angular mode is
\[
 \Lop_{\delta,j}
 =-\frac{d^2}{dr^2}-\frac1r\frac d{dr}
   +\frac{j^2}{r^2}+1-pU_\delta^{p-1}.
\]
Set
\begin{equation}\label{eq:rho-scaling}
 \varepsilon=\frac1K,\qquad
 \rho=\frac{K\delta}{\lambda},\qquad
 \delta=\lambda\varepsilon\rho,\qquad
 w_{\varepsilon,\rho}(t)=1+\varepsilon\rho t.
\end{equation}
We use $t=(r-\lambda)/\delta$ and write $w=w_{\varepsilon,\rho}$.
Then $r=\lambda w$, $d/dr=\delta^{-1}d/dt$, and
$\delta^2U_\delta^{p-1}=V_\delta^{p-1}$. For $j=K$, we also have
$\delta^2K^2/r^2=\rho^2/w^2$. Thus multiplication of the
original eigenvalue equation by $\delta^2$ gives the operator
\begin{equation}\label{eq:scaled-K}
 A_{\varepsilon,\rho}h
 =-\frac1w(wh')'
  +\left(\frac{\rho^2}{w^2}
        +\lambda^2\varepsilon^2\rho^2
        -pV_{\lambda\varepsilon\rho}^{p-1}\right)h,
\end{equation}
with mixed boundary conditions
\begin{equation}\label{eq:scaled-mixed-bc}
 h(0)=0,\qquad h'(1)
       +\frac{\varepsilon\rho}{1+\varepsilon\rho}h(1)=0.
\end{equation}
Indeed, the original Robin condition is
$\delta^{-1}h'(1)+R^{-1}h(1)=0$, and
$\delta/R=\varepsilon\rho/(1+\varepsilon\rho)$.
Primes in \eqref{eq:scaled-K}--\eqref{eq:scaled-mixed-bc} denote
derivatives with respect to $t$.

We use the fixed form domain
$E=\{h\in H^1(0,1):h(0)=0\}$ and the weighted norm form
\[
 n_{\varepsilon,\rho}(h)=\int_0^1wh^2\dd t.
\]
The corresponding quadratic form is
\begin{align}
 q_{\varepsilon,\rho}(h)
 ={}&\int_0^1 w
 \left\{|h'|^2+
 \left(\frac{\rho^2}{w^2}
 +\lambda^2\varepsilon^2\rho^2
 -pV_{\lambda\varepsilon\rho}^{p-1}\right)h^2\right\}\dd t
 \notag\\
 &+\varepsilon\rho\,h(1)^2.
 \label{eq:scaled-form}
\end{align}
To check the endpoint term, integration by parts gives
\[
 \int_0^1wh\,A_{\varepsilon,\rho}h\dd t
 =\int_0^1w\left\{|h'|^2+
 \left(\frac{\rho^2}{w^2}+\lambda^2\varepsilon^2\rho^2
       -pV_{\lambda\varepsilon\rho}^{p-1}\right)h^2\right\}\dd t
       -w(1)h(1)h'(1).
\]
The left endpoint contributes zero because $h(0)=0$.
The Robin condition makes the last term equal to
$\varepsilon\rho h(1)^2$, since its coefficient multiplied by
$w(1)$ is $\varepsilon\rho$. The form domain therefore does not
vary when the Robin coefficient varies.

Denote by $e_\varepsilon(\rho)$ the principal eigenvalue, defined
by the minimum of $q_{\varepsilon,\rho}/n_{\varepsilon,\rho}$
on $E\setminus\{0\}$. The interval is bounded, the weights are
strictly positive, and the potential is bounded. The compact
embedding $E\hookrightarrow L^2(0,1)$ gives a minimizer.
The associated separated Sturm--Liouville problem has a simple
principal eigenvalue and a positive eigenfunction.
At $\varepsilon=0$, the limiting form is
\[
 q_{0,\rho}(h)=
 \int_0^1\bigl(|h'|^2+(\rho^2-pV_p^{p-1})h^2\bigr)\dd t.
\]
Its natural endpoint condition is $h'(1)=0$. Addition of the
constant potential $\rho^2$ shifts every mixed eigenvalue by
$\rho^2$, so
\begin{equation}\label{eq:limit-eigenvalue}
 e_0(\rho)=\nu_{\mathrm{DN}}+\rho^2=\rho^2-\kappa_p^2.
\end{equation}

\subsubsection{Uniform \texorpdfstring{$C^1$}{C1} convergence of the principal eigenvalue}
Fix a compact interval $J\Subset(0,\infty)$ whose interior
contains $\kappa_p$. All estimates in this subsubsection are
uniform for $\rho\in J$.
We next determine the limiting spectral crossing.

\begin{lemma}\label{lem:C1-crossing}
As $\varepsilon\downarrow0$,
\begin{equation}\label{eq:C1-eigen}
 e_\varepsilon(\rho)\longrightarrow\rho^2-\kappa_p^2,\qquad
 e_\varepsilon'(\rho)\longrightarrow2\rho
 \quad\text{uniformly on }J.
\end{equation}
\end{lemma}

\begin{proof}
Since $\delta=\lambda\varepsilon\rho$ tends uniformly to zero,
Lemma~\ref{lem:potential-derivative} gives uniform convergence
of the potentials. Also $w\to1$ uniformly and
$\varepsilon\rho\to0$. The trace estimate
$|h(1)|^2\le C\|h\|_{H^1}^2$ therefore gives numbers
$\omega_\varepsilon\downarrow0$ such that
\[
 |q_{\varepsilon,\rho}(h)-q_{0,\rho}(h)|
       \le\omega_\varepsilon\|h\|_{H^1}^2,\qquad
 |n_{\varepsilon,\rho}(h)-\|h\|_2^2|
       \le C\varepsilon\|h\|_2^2.
\]
The bounded potentials and the positive weights also imply
\[
 q_{\varepsilon,\rho}(h)
       \ge c\|h'\|_2^2-C\|h\|_2^2,
 \qquad c>0,
\]
with constants independent of small $\varepsilon$ and of $\rho$.
Testing the Rayleigh quotient with a fixed normalized principal
eigenfunction of the limiting mixed operator gives
$e_\varepsilon(\rho)\le e_0(\rho)+o(1)$ uniformly.
In particular, the principal eigenvalues are uniformly bounded.
Normalize their positive eigenfunctions
$\phi_{\varepsilon,\rho}$ by
$n_{\varepsilon,\rho}(\phi_{\varepsilon,\rho})=1$.
The lower form bound then gives a uniform $H^1$ bound for these
eigenfunctions. Evaluating the form comparison on them and using
$q_{0,\rho}(h)\ge e_0(\rho)\|h\|_2^2$ proves the reverse
inequality up to $o(1)$. This establishes the first convergence
in \eqref{eq:C1-eigen}.

We next justify differentiation of the eigenvalue. For fixed
$\varepsilon>0$, write
$V_{\lambda\varepsilon\rho}(t)=(1-t)z_\rho(t)$.
The quotient is strictly positive on $[0,1]$ and depends smoothly
on $\rho$ in $C^{1,\alpha}$. Hence
\[
 V_{\lambda\varepsilon\rho}^{p-1}
       =(1-t)^{p-1}z_\rho^{p-1}
\]
depends smoothly on $\rho$ in $C^{0,\alpha}$. This factorization
also covers noninteger $p$ and $1<p<2$, when differentiation of
the uncombined power would produce a singular factor at $t=1$.
The forms are thus smooth on the fixed domain $E$. Simplicity of
the principal eigenvalue gives locally differentiable normalized
eigenfunctions and a differentiable eigenvalue.

For clarity, let $q_{\varepsilon,\rho}(h,k)$ and
$n_{\varepsilon,\rho}(h,k)$ denote the bilinear forms. The weak
eigenvalue equation is
\[
 q_{\varepsilon,\rho}(\phi,k)
       =e_\varepsilon(\rho)n_{\varepsilon,\rho}(\phi,k),
       \qquad k\in E.
\]
Differentiating $q(\phi,\phi)=e_\varepsilon n(\phi,\phi)$
and using this equation with $k=\partial_\rho\phi$ cancels
the two terms involving $\partial_\rho\phi$. Since
$n(\phi,\phi)=1$, we obtain
\begin{equation}\label{eq:HF}
 e_\varepsilon'(\rho)
 =\partial_\rho q_{\varepsilon,\rho}
      (\phi_{\varepsilon,\rho})
 -e_\varepsilon(\rho)
       \int_0^1(\partial_\rho w)\phi_{\varepsilon,\rho}^2\dd t.
\end{equation}
In particular, this identity includes the change in the weighted
normalization. It also includes the Robin variation through the
endpoint term in $q$.

We now estimate each coefficient derivative explicitly. Put
$P_{\varepsilon,\rho}=V_{\lambda\varepsilon\rho}^{p-1}$.
Since $\partial_\rho w=\varepsilon t$, differentiation of
\eqref{eq:scaled-form} at a fixed $h$ gives
\begin{align*}
 \partial_\rho q_{\varepsilon,\rho}(h)
 ={}&\int_0^1\varepsilon t|h'|^2\dd t
        +\varepsilon h(1)^2\\
 &+\int_0^1\left\{
       \frac{2\rho}{w}-\frac{\varepsilon\rho^2t}{w^2}
       +\lambda^2\varepsilon^2(2\rho w+\varepsilon\rho^2t)
       -p\varepsilon tP_{\varepsilon,\rho}
       -pw\partial_\rho P_{\varepsilon,\rho}
                    \right\}h^2\dd t.
\end{align*}
The derivative of the profile is taken along the radial branch:
\[
 \partial_\rho P_{\varepsilon,\rho}
 =\lambda\varepsilon\partial_\delta(V_\delta^{p-1})
 =\frac1\rho\delta\partial_\delta(V_\delta^{p-1}).
\]
Because $J$ is bounded away from zero,
Lemma~\ref{lem:potential-derivative} shows that its $L^\infty$
norm tends uniformly to zero. The angular coefficient in braces
therefore tends uniformly to $2\rho$, and all other displayed
terms are $o(1)\|h\|_{H^1}^2$. Thus
\[
 \partial_\rho q_{\varepsilon,\rho}(h)
       =2\rho\|h\|_2^2+o(1)\|h\|_{H^1}^2.
\]
The normalized eigenfunctions have uniformly bounded $H^1$ norms,
and $w\to1$ implies $\|\phi_{\varepsilon,\rho}\|_2^2\to1$
uniformly. Finally, the last term in \eqref{eq:HF} is $O(\varepsilon)$.
Substitution proves the derivative convergence and completes the
proof.
\end{proof}

Choose $0<b<\kappa_p/2$ and take
$J=[\kappa_p-b,\kappa_p+b]$. The limiting eigenvalue is negative
at the left endpoint and positive at the right endpoint. Its
derivative $2\rho$ is strictly positive on $J$.
Lemma~\ref{lem:C1-crossing} therefore implies that, for all
sufficiently small $\varepsilon$, $e_\varepsilon$ has the same
endpoint signs and is strictly increasing throughout $J$.
The intermediate value theorem gives a unique zero
$\rho_\varepsilon\in J$, and uniform convergence gives
\begin{equation}\label{eq:rho-root}
 \rho_\varepsilon\longrightarrow\kappa_p,\qquad
 e_\varepsilon'(\rho_\varepsilon)>0.
\end{equation}
Taking $\varepsilon=1/K$, we define
\begin{equation}\label{eq:critical-radius}
 \delta_K=\frac{\lambda\rho_{1/K}}K,\qquad
 R_K=\lambda+\delta_K.
\end{equation}
The original mixed eigenvalue is
$\zeta_K(R)=\delta^{-2}e_{1/K}(\rho)$, so these radii give
its required zeros. Equation~\eqref{eq:RK-asymptotic} follows
immediately. The uniqueness just proved concerns the chosen radial
family and the stated interval of the thin-annulus scaling parameter.

\subsubsection{Dirichlet invertibility in the symmetry class}
Let $e^D_\varepsilon(\rho)$ be the principal eigenvalue of
the same operator with $h(0)=h(1)=0$. Its form domain is
$H_0^1(0,1)$, so the endpoint term vanishes. The preceding form
comparison on this smaller domain gives
\begin{equation}\label{eq:D-K-limit}
 e^D_\varepsilon(\rho)\longrightarrow\nu_D+\rho^2.
\end{equation}
In particular, at $\rho=\rho_\varepsilon$ the limiting value is
$\nu_D+\kappa_p^2=g_p>0$. This is the strict gap from
Lemma~\ref{lem:gap}: the mixed boundary problem can have a zero
eigenvalue while the Dirichlet problem in the same angular mode
remains strictly positive.
The following proposition establishes the Dirichlet gap at the nonlinear critical radius.

\begin{proposition}
\label{prop:Dirichlet-gap}
For every sufficiently large $K$, the Dirichlet realization of
$-\Delta+1-pU_K^{p-1}$ on $A_{\lambda,R_K}$ is invertible in
the space of $D_K$-invariant functions. It has exactly one negative
eigenvalue in this space, belonging to the radial mode. The quadratic
form is positive in every permitted nonradial Dirichlet mode.
\end{proposition}

\begin{proof}
For the radial mode, the scaled operator is
\eqref{eq:radial-Jacobi-scaled}. Its first two Dirichlet eigenvalues
converge to the first two eigenvalues of $A_D$. By
Lemma~\ref{lem:pulse-D}, the first limiting eigenvalue is negative
and the second is positive. Hence, for small $\delta_K$, the radial
operator has one negative eigenvalue and a strictly positive second
eigenvalue. All its remaining eigenvalues are at least the second
one, so there is no radial zero eigenvalue.

For $j=K$, \eqref{eq:D-K-limit} and
$\rho_{1/K}\to\kappa_p$ give
$e^D_{1/K}(\rho_{1/K})\ge g_p/2$ for sufficiently large $K$.
For $j=nK$, $n\ge2$, the scaled operator has, in addition, the
potential
\[
 \frac{(n^2-1)\rho_{1/K}^2}
      {(1+\rho_{1/K}t/K)^2}\ge0.
\]
Its Rayleigh quotient is therefore at least that of the $K$-mode
for every test function. This gives the same lower bound $g_p/2$
for the principal scaled Dirichlet eigenvalue in every permitted
nonradial sector. Every higher eigenvalue in such a sector is also
positive. Multiplication by $\delta_K^{-2}$ returns these bounds
to the original variables without changing their signs.

A $D_K$-invariant function has only the radial mode and the cosine
modes $nK$, $n\ge1$. These sectors are orthogonal both for the
$L^2$ norm and for the Dirichlet quadratic form, because the
coefficients are radial. The preceding bounds thus give a trivial
kernel and exactly one negative eigenvalue on the symmetry class.
The potential $1-pU_K^{p-1}$ belongs to $C^{0,\alpha}$ by the
outer-boundary factorization of $U_K$ and the restriction on
$\alpha$. The Dirichlet Schauder realization is Fredholm of index
zero. Its trivial kernel gives surjectivity, and the Schauder
estimate and bounded inverse theorem give a bounded isomorphism
on the stated Banach spaces.
\end{proof}

The proposition concerns the symmetry class used in the bifurcation
argument. This class excludes the lower angular modes
$1,\ldots,K-1$. The state equation requires Dirichlet nonresonance
in this class; the argument does not assert full-space nonresonance.

\subsubsection{Isolation of the boundary kernel}
At $R_K$, the principal mixed eigenvalue of the $K$-operator is
zero and is simple. Its positive eigenfunction has a nonzero outer
value: a zero outer value together with the Robin condition would
give zero Cauchy data at $R_K$, contradicting ODE uniqueness.
We can therefore normalize this eigenfunction to have value $1$
at $R_K$. By the Dirichlet uniqueness just proved, it is precisely
the function $H_K$ in \eqref{eq:Hj}. Its Robin condition gives
$T_K(R_K)=0$.

For every $n\ge2$, \eqref{eq:mode-order} gives
\[
 \zeta_{nK}(R_K)\ge\frac{(n^2-1)K^2}{R_K^2}>0.
\]
Since the principal eigenvalue is already positive in each of these
sectors, zero cannot be any of its higher mixed eigenvalues.
Equivalently, $T_{nK}(R_K)\ne0$. The boundary variable is
mean-zero, so it has no radial mode. Lemma~\ref{lem:power-Fredholm}
therefore gives the one-dimensional boundary kernel and its stated
codimension-one range. It remains to transfer the positive
derivative in \eqref{eq:rho-root} to the boundary symbol.

\subsubsection{Transversality of the boundary operator}
Fix a sufficiently large $K$ and put $\varepsilon=1/K$.
For $\rho$ near $\rho_{1/K}$ and $z$ near zero, let
$h_{z,\rho}$ solve
\[
 (A_{1/K,\rho}-z)h_{z,\rho}=0,\qquad
 h_{z,\rho}(0)=0,\quad h_{z,\rho}(1)=1.
\]
The positive Dirichlet gap persists in a neighborhood of the
critical parameters. Subtracting a fixed lift of the endpoint
values reduces this problem to an invertible zero-Dirichlet
problem, which gives existence, uniqueness, and smooth dependence
on $(z,\rho)$. Define
\begin{equation}\label{eq:Bzrho}
 B(z,\rho)=h_{z,\rho}'(1)
        +\frac{\rho/K}{1+\rho/K}.
\end{equation}
The radial solution and the coefficients of $A_{1/K,\rho}$ vary
with $\rho$ in this definition. Returning to the original radial
variable in \eqref{eq:Hj}, we obtain
\begin{equation}\label{eq:T-B}
 T_K(R)=\frac1\delta B(0,\rho),\qquad
 \rho=\frac{K(R-\lambda)}{\lambda}.
\end{equation}
Indeed, $H_K'(R)=\delta^{-1}h_{0,\rho}'(1)$ and
$\delta/R=(\rho/K)/(1+\rho/K)$.
The principal mixed eigenfunction has nonzero right endpoint value,
so its normalization to $1$ also gives
$B(e_{1/K}(\rho),\rho)=0$ throughout a sufficiently small
neighborhood of the critical point.

We first compute the sign of $B_z$. Write $h=h_{z,\rho}$
and $\dot h=\partial_z h_{z,\rho}$. Differentiating the equation
and the fixed endpoint data gives
\[
 (A_{1/K,\rho}-z)\dot h=h,\qquad
 \dot h(0)=\dot h(1)=0.
\]
Multiply this equation by $wh$, multiply the equation for $h$
by $w\dot h$, and subtract. The potential and spectral terms
cancel. The remaining Green identity is
\begin{align*}
 \int_0^1wh^2\dd t
 &=\int_0^1w\bigl(h(A-z)\dot h-\dot h(A-z)h\bigr)\dd t\\
 &=\bigl[w(\dot h h'-h\dot h')\bigr]_0^1
 =-w(1)\dot h'(1).
\end{align*}
Both contributions at $t=0$ vanish because $h(0)=\dot h(0)=0$.
At $t=1$, the term containing $\dot h(1)$ vanishes, while
$h(1)=1$. The Robin coefficient is independent of $z$, hence
\begin{equation}\label{eq:Bz-sign}
 B_z(z,\rho)
 =-\frac{\int_0^1wh^2\dd t}{w(1)}<0.
\end{equation}
This calculation in fact applies to every $(z,\rho)$ in the
above Dirichlet-invertible neighborhood.

We next differentiate $B(e_{1/K}(\rho),\rho)=0$. At the
critical point $e_{1/K}(\rho_{1/K})=0$, the chain rule gives
\begin{equation}\label{eq:B-rho}
 B_\rho(0,\rho_{1/K})
 =-B_z(0,\rho_{1/K})e_{1/K}'(\rho_{1/K})>0.
\end{equation}
Here $B_\rho$ is the partial derivative with $z$ fixed. It contains
both the derivative of $h_{z,\rho}'(1)$ and the derivative
$K^{-1}(1+\rho/K)^{-2}$ of the Robin coefficient.
One can also see these terms directly from Green's identity.
Writing $\beta=(\rho/K)/(1+\rho/K)$, differentiation with
respect to $\rho$ at fixed $z$ gives
\[
 (A-z)h_\rho=-A_\rho h,\qquad h_\rho(0)=h_\rho(1)=0,
 \qquad
 w(1)h_\rho'(1)=\int_0^1wh A_\rho h\dd t.
\]
At a mixed eigenvalue, $h'(1)=-\beta$. Differentiation of the
form at a fixed $h$ and integration by parts therefore give
\[
 \int_0^1wh A_\rho h\dd t
 =\partial_\rho q_{1/K,\rho}(h)
     -z\partial_\rho n_{1/K,\rho}(h)-w(1)\beta_\rho.
\]
The term $-w_\rho(1)h(1)h'(1)$ from integration by parts
cancels the $w_\rho(1)\beta$ part of the differentiated endpoint
term. Adding $w(1)\beta_\rho$ thus yields
\[
 B_\rho(z,\rho)
 =\frac{\partial_\rho q_{1/K,\rho}(h)
       -z\partial_\rho n_{1/K,\rho}(h)}{w(1)}
 =\frac{n_{1/K,\rho}(h)}{w(1)}e_{1/K}'(\rho),
 \qquad z=e_{1/K}(\rho).
\]
The last equality is \eqref{eq:HF} with an unnormalized
eigenfunction. It agrees with \eqref{eq:B-rho} and explicitly
accounts for the weight, profile, and Robin variations.

Since $d\rho/dR=K/\lambda$ and $d\delta/dR=1$,
differentiation of \eqref{eq:T-B} gives
\[
 T_K'(R)=-\frac1{\delta^2}B(0,\rho)
                 +\frac{K}{\lambda\delta}B_\rho(0,\rho).
\]
The first term vanishes at $R_K$. Combining
\eqref{eq:Bz-sign} and \eqref{eq:B-rho} gives
\begin{equation}\label{eq:T-radius-positive}
 T_K'(R_K)
 =\frac{K}{\lambda\delta_K}
   \frac{\int_0^1w h_{0,\rho_{1/K}}^2\dd t}{w(1)}
       e_{1/K}'(\rho_{1/K})>0.
\end{equation}
Finally, \eqref{eq:multiplier} gives
$m_K(R)=-2d_{R-\lambda}^2T_K(R)$. Its derivative contains
the terms $-4d_{R-\lambda}(\partial_Rd_{R-\lambda})T_K$
and $-2d_{R-\lambda}^2T_K'$. The first term vanishes at the
critical radius because $T_K(R_K)=0$. Thus
\begin{equation}\label{eq:nonlinear-transverse}
 m_K'(R_K)=-2d_{\delta_K}^2T_K'(R_K)<0.
\end{equation}
To express the Crandall--Rabinowitz condition directly, put
$\phi_K(\theta)=\cos(K\theta)$ and
$\ell_K(g)=\int_0^{2\pi}g(\theta)\phi_K(\theta)\dd\theta$.
Lemma~\ref{lem:power-Fredholm} gives
$\Ran D_\eta F_{p,K}(R_K,0)=\Ker\ell_K$. Differentiating
\eqref{eq:multiplier} along the radial branch now yields
\[
 \ell_K\bigl(D_RD_\eta F_{p,K}(R_K,0)[\phi_K]\bigr)
 =m_K'(R_K)\int_0^{2\pi}\cos^2(K\theta)\dd\theta
 =\pi m_K'(R_K)<0.
\]
Consequently $D_RD_\eta F_{p,K}(R_K,0)[\phi_K]$ does not
belong to the range. This proves transversality with respect to
the outer radius, including the changes in the endpoint, radial
solution, weight, and Robin coefficient.

\subsection{Application of the Crandall--Rabinowitz theorem}
We first fix the threshold for the modes independently of the H\"older
exponent. Carry out the radial construction and spectral estimates
with $\alpha_*=\frac12\min\{1,p-1\}$. The thresholds then depend
only on $\lambda,p$. For fixed $\delta>0$, the radial equation
implies $V_\delta\in C^3([0,1])$, since $s\mapsto s^p$ is
$C^1$ on $[0,\infty)$. Therefore $V_\delta\in C^{2,\alpha}$
for every $\alpha<1$. For each exponent in \eqref{eq:alpha},
Dirichlet nondegeneracy, the Schauder estimates, and
Lemma~\ref{lem:power-map} give the same smooth local radial family
in $C^{2,\alpha}$. Uniqueness in the weaker of the two spaces
identifies the families. We may thus choose
$K_0=K_0(\lambda,p)\ge2$ uniformly for these exponents, although
the local neighborhoods may depend on $\alpha$.

Fix $K\ge K_0$, set $R_*=R_K$ in
\eqref{eq:power-pullback}, and let $U_K=U_{\delta_K}$. By
Proposition~\ref{prop:Dirichlet-gap}, the power-law state
linearization is an isomorphism $Z_{p,K}\to W_{p,K}$.
Proposition~\ref{prop:power-state} therefore supplies the smooth
positive state map used in \eqref{eq:power-Bernoulli}.
All hypotheses of the Crandall--Rabinowitz theorem \cite{CR} can now
be verified for this concrete map:
\begin{enumerate}
\item $F_{p,K}$ is $C^\infty$ from a neighborhood of $(R_K,0)$
      in $\R\times X_K$ into $Y_K$, and $F_{p,K}(R,0)=0$.
\item By Lemma~\ref{lem:power-Fredholm} and the isolation of the
      principal mixed mode,
      $L_K=D_\eta F_{p,K}(R_K,0)$ is Fredholm of index zero, with
      $\Ker L_K=\Span\{\phi_K\}$, where
      $\phi_K(\theta)=\cos(K\theta)$, and
      $\Ran L_K=\{g\in Y_K:\int_0^{2\pi}g\phi_K=0\}$.
\item By \eqref{eq:nonlinear-transverse},
      \[
      D_RD_\eta F_{p,K}(R_K,0)\phi_K
      =m_K'(R_K)\phi_K\notin\Ran L_K.
      \]
\end{enumerate}
The theorem gives a $C^\infty$ curve
\[
 R=R(s),\qquad \eta(s)=s\phi_K+s z(s),\qquad z(0)=0,
\]
with $z(s)$ in the fixed complement
\[
 X_K^\perp=\left\{\eta\in X_K:
      \int_0^{2\pi}\eta(\theta)\phi_K(\theta)\dd\theta=0\right\}.
\]
We use the Fourier coefficient
$s=\pi^{-1}\int_0^{2\pi}\eta(s,\theta)\phi_K(\theta)\dd\theta$
as the curve parameter. The half-period rotation
$\mathcal R\eta(\theta)=\eta(\theta+\pi/K)$ preserves both
the equation and this complement, while reversing this coefficient.
Local uniqueness of the bifurcation curve therefore yields
\begin{equation}\label{eq:branch-parity}
 R(-s)=R(s),\qquad \eta(-s)=\mathcal R\eta(s),\qquad Q(-s)=Q(s).
\end{equation}
Here $Q(s)$ is defined by \eqref{eq:power-recover-Q}. Consequently,
\[
 R(s)=R_K+O(s^2),\qquad
 \eta(s)=s\cos(K\theta)+O(s^2),\qquad Q(s)>0
\]
for sufficiently small $|s|$.

The positivity in Proposition~\ref{prop:power-state} carries over
to the physical solutions. Moreover,
$U_K'(\lambda)>0$, and the coordinate map fixes an inner collar.
Continuity in $C^1$ gives
$\partial_r u_s(\lambda,\theta)>0$ for small $|s|$.
Since the inner value is $1$, each solution exceeds $1$ in an
inner collar and thus lies outside the ordered class.
For $s\ne0$, the $K$-th Fourier coefficient of $\eta(s)$ is
nonzero. If the outer boundary were a circle, its rotational
invariance by $2\pi/K$, with $K\ge2$, would force its center
to be the origin. Its graph would then be constant, and the zero-mean
condition would give $\eta=0$, a contradiction. Thus the boundary
is noncircular.

We next verify the local completeness statement, including the state
symmetry. Suppose that a nearby positive solution of
\eqref{eq:power} has boundary perturbation $\eta\in X_K$.
For $g\in D_K$, the states $u$ and $u\circ g$ are defined on
the same domain and have the same Dirichlet data. Positivity and the
outer Bernoulli condition imply
$\partial_\nu u=-\sqrt Q$, so they have the same outer normal
derivative as well. Their difference $w$ satisfies
\[
 \Delta w=c_g(x)w,\qquad
 c_g(x)=1-p\int_0^1
  \bigl(tu(x)+(1-t)u(gx)\bigr)^{p-1}\dd t,
\]
where $c_g$ is bounded, and has zero Dirichlet and Neumann traces on
the outer boundary. Extending $w$ by zero through a small outer
collar gives a weak solution with bounded potential, since the two
traces produce no boundary distribution. Unique continuation for
elliptic equations with bounded zeroth-order coefficients
\cite{Aronszajn} implies $w=0$. Therefore the state is $D_K$ invariant.
The uniqueness in Proposition~\ref{prop:power-state} identifies its
pullback with $v_{R,\eta}$, and its boundary satisfies
$F_{p,K}(R,\eta)=0$. The local completeness part of the
Crandall--Rabinowitz theorem now puts every nearby nonradial solution
on the curve just constructed. The neighborhood is taken in both the
boundary variables and the pulled-back $C^{2,\alpha}$ state.

Finally, \eqref{eq:rho-root}--\eqref{eq:critical-radius} give
\[
 R_K=\lambda+\frac{\lambda\kappa_p}{K}+o(K^{-1}).
\]
This proves Theorem~\ref{thm:power}.

\subsection{Further properties of the power-law solutions}
\label{subsec:power-properties}

\begin{proof}[Proof of Proposition~\ref{prop:first}]
Set $d=U'(R)$. Since the solution is positive and $U(R)=0$,
we have $d\le0$. If $d=0$, uniqueness for the radial Cauchy
problem at $r=R$ would give the zero solution, contradicting
$U(\lambda)=1$. Therefore $d<0$.

Differentiating the radial equation, we find that $U'$ satisfies
the first angular Jacobi equation
\begin{equation}\label{eq:translation-Jacobi}
 -(U')''-\frac1r(U')'
 +\left(\frac1{r^2}+1-pU^{p-1}\right)U'=0.
\end{equation}
The identity $U''(R)=-d/R$ shows that $U'$ also satisfies the
outer Robin condition. Suppose that a nonzero mixed solution $h$
exists. Its outer value must be nonzero, since otherwise the Robin
condition would give zero Cauchy data. Normalize it by $h(R)=1$.
Its outer Cauchy data then coincide with those of $U'/d$, so ODE
uniqueness gives $h=U'/d$ on the full interval. The inner condition
$h(\lambda)=0$ consequently implies $U'(\lambda)=0$.

Since $1-1^p=0$, the constant $u=1$ is an equilibrium solution
of the radial equation. Uniqueness for the initial data
$U(\lambda)=1$ and $U'(\lambda)=0$ therefore gives
$U\equiv1$, contradicting $U(R)=0$. This proves the proposition.
\end{proof}

\bigskip

\begin{center}
\textbf{Acknowledgement}
\end{center}

The research of Wang was supported by supported by National
Key R$\&$D Program of China 2022YFA1005601 and National Natural
Science Foundation of China 12371114. The research of Xu was supported by National Natural Science Foundation of China (No. 12401134). The research of Zhang was
supported by National Natural Science Foundation of China (No. 12301133) and the Postdoctoral Science
Foundation of China (No. 2024T170353).

\end{document}